\documentclass[onefignum,onetabnum,pagebackref]{siamart220329}

\usepackage{lipsum}
\usepackage{amsfonts}
\usepackage{graphicx}
\usepackage{epstopdf}
\usepackage{tikz}
\usetikzlibrary{positioning}
\ifpdf
  \DeclareGraphicsExtensions{.eps,.pdf,.png,.jpg}
\else
  \DeclareGraphicsExtensions{.eps}
\fi

\usepackage{enumitem}
\setlist[enumerate]{leftmargin=.5in}
\setlist[itemize]{leftmargin=.5in}

\usepackage{amsmath, amssymb, mathtools, stmaryrd}
\usepackage{amsfonts,mathrsfs, mleftright,booktabs,soul}
\usepackage{subcaption}
\usepackage{xspace}
\usepackage{algorithm, algpseudocode, algorithmicx}
\usepackage{backref}

\newsiamthm{fact}{Fact}
\newsiamthm{question}{Question}

\headers{High-accuracy randomized nullspace}{E. N. Epperly, T. Park, and Y. Nakatsukasa}

\title{Fast, high-accuracy, randomized nullspace computations for tall matrices\thanks{Date: \today}\funding{The authors acknowledge support from the Miller Institute for Basic Research in Science, University
of California Berkeley (ENE), the RandESC project funded by the Swiss Platform for Advanced Scientific Computing (TP), and EPSRC grants EP/Y010086/1 and EP/Y030990/1 (YN).
ENE also received support, under aegis of Joel Tropp, by ONR grant N00014-24-1-2223 and the Caltech
Carver Mead New Adventures Fund.}}

\author{Ethan N. Epperly\thanks{Department of Mathematics, University of California Berkeley, Berkeley, CA 94720, USA (\url{https://ethanepperly.com}, \email{eepperly@berkeley.edu})}
\and Taejun Park\thanks{Institute of Mathematics, EPF Lausanne, 1015 Lausanne, Switzerland, (\email{taejun.park@epfl.ch}).} \and Yuji Nakatsukasa\thanks{Mathematical Institute, University of Oxford, Oxford, OX2 6GG, UK, (\email{nakatsukasa@maths.ox.ac.uk}).} 
}

\usepackage{amsopn}

\renewcommand*{\backref}[1]{}
\renewcommand*{\backrefalt}[4]{%
	\ifcase #1 %
	(No citations.)
	\or
	(Cited p.~#2.)
	\else
	(Cited pp.~#2.)
	\fi
}

\usepackage{xcolor}
\definecolor{Mulberry}{HTML}{C54B8C}

\algrenewcommand\algorithmiccomment[1]{\hfill\textcolor{Mulberry}{$\triangleright$ #1}}

\newcommand{\real}{\mathbb{R}}
\newcommand{\complex}{\mathbb{C}}

\renewcommand{\left}{\mleft}
  \renewcommand{\right}{\mright}

\DeclareMathOperator{\diag}{diag}
\newcommand{\mat}[1]{\boldsymbol{#1}}
\renewcommand{\vec}[1]{\boldsymbol{#1}}

\newcommand{\norm}[1]{\left\| #1 \right\|}
\DeclareMathOperator{\spn}{span}

\newcommand{\QR}{\textsf{QR}\xspace}
\newcommand{\evec}{\mathbf{e}}

\newcommand{\Id}{\mathbf{I}}

\newcommand{\order}{\mathcal{O}}

\DeclareMathOperator{\poly}{poly}

\DeclareMathOperator*{\argmin}{argmin}

\renewcommand{\Re}{\mathrm{Re}}

\renewcommand{\hat}[1]{\widehat{#1}}
\renewcommand{\tilde}[1]{\widetilde{#1}}

\newcommand{\gap}{\mathrm{gap}}

\usepackage[sort]{cite}

\usepackage[normalem]{ulem}

\newcommand{\ignore}[1]{}

\begin{document}
	
\maketitle
	
\begin{abstract}    In this paper, we develop RLOBPCG, an efficient method for computing a small number of singular triplets corresponding to the smallest singular values of large, tall matrices. The algorithm combines randomized preconditioner from the sketch-and-precondition techniques with the LOBPCG eigensolver: a small sketch is used to construct a high-quality preconditioner, and LOBPCG is run on the Gram matrix to refine the singular vector. Under the standard subspace embedding assumption and a modest singular value gap between the two smallest singular values, we prove that RLOBPCG converges geometrically to the minimum singular vector. In numerical experiments, RLOBPCG achieves near-optimal accuracy on matrices with up to $10^6$ rows, outperforming classical LOBPCG and Lanczos methods by a speedup of up to $12\times$ and maintaining robustness when other iterative methods fail to converge. 
\end{abstract}
	
	
	
\section{Introduction} 
The singular vectors of a matrix $\mat{A}\in\complex^{m\times n}$ are fundamental objects in matrix computations and data analysis.
Much work has been dedicated to efficiently computing the \emph{dominant} singular vectors, those associated with the largest singular values.
The dominant singular values comprise the optimal low-rank approximation to $\mat{A}$, and computing them is necessary for applications in genetic data analysis \cite{AI14}, model reduction \cite{AK19a}, compression of large data \cite{TYUC19}, and many other areas.

Computation of the \emph{minimum singular vectors}, those with the smallest singular values, has received significantly less attention.
In the case where $\mat{A}$ is tall and rank-deficient, the minimum singular vectors span the \emph{nullspace} of $\mat{A}$.
When $\mat{A}$ is not rank-deficient, the minimum singular vectors can be interpreted as spanning the \emph{near-nullspace} of $\mat{A}$, the vectors that are amplified by $\mat{A}$ the least.
Minimum singular vector computations are used in applications spanning total least squares \cite{GV80}, rational approximation \cite{nakatsukasa2018aaa,nakatsukasa2020algorithm}, frequency estimation \cite{Sch86}, computational topology \cite{KS24}, and multivariate polynomial root-finding \cite{BDD14,GWCD24}.

This paper focuses on computing a small number, say 1 to 10, of minimum singular vectors of a large, very tall matrix $\mat{A} \in \complex^{m\times n}$, for which $m\gg n\gg 1$.
Our approach combines the \emph{randomized preconditioning} approach of Rokhlin \& Tygert \cite{RT08} (popularized as \emph{Blendenpik} by \cite{AMT10}) with a preconditioned eigensolver such as LOBPCG \cite{Kny01}.
This combination is simple and natural, but to the best of our knowledge, has not yet been proposed, analyzed, and tested.
We use the term RLOBPCG to refer to the LOBPCG method with randomized preconditioning. In this paper, we study the use of randomized preconditioning to accelerate the LOBPCG algorithm.
The choice of LOBPCG in this proposal is not essential, and one can also use randomized preconditioning to preconditioner other eigensolvers or singular value solvers such as preconditioned inverse iteration (PINVIT, \cite{KN03}), preconditioned steepest descent \cite{Ney12}, Jacobi--Davidson SVD (JDSVD, \cite{Hoc01}), PHSVDS \cite{WS15}, and Golub--Kahan Davidson \cite{GSR19}.
Our choice of LOBPCG in this paper was motivated by preliminary experiments, where we found RLOBPCG to perform much better than PINVIT and preconditioned descent with randomized preconditioning.
The investigation of randomized preconditioning with dedicated singular value solvers \cite{Hoc01,WS15,GSR19} is a natural subject for future work.

To compute a single minimum singular vector, RLOBPCG proceeds as follows:
\begin{enumerate}
    \item \textbf{Generate embedding.} Draw a random embedding matrix (i.e., a ``sketching matrix'') $\mat{S}\in\complex^{d\times m}$ with $d = \order(n)$ or $d = \order(n\log n)$ rows.
    \item \textbf{Sketch and factorize.} Compute the sketch $\mat{S}\mat{A}$ and form its (economy-size) decomposition $\mat{S}\mat{A} = \mat{U}\mat{\Sigma} \mat{V}^*$. 
    (${}^*$ denotes the conjugate transpose.)
    \item \textbf{Define preconditioner and initialization.} Set $\mat{P} \coloneqq \mat{V}\mat{\Sigma}^{-1}$ and $\vec{v}_0 = \vec{v}_{-1} \coloneqq \mat{P}\evec_n / \norm{\mat{P}\evec_n}$, where $\evec_n$ is the $n$th standard basis vector.\footnote{More traditionally, randomized preconditioning schemes typically use a \QR decomposition of the sketch $\mat{S}\mat{A}$ rather than an SVD.
    We use an SVD here because it immediately yields the good initializations $\vec{v}_0 = \vec{v}_{-1}$.
    We will also use the singular value matrix $\mat{\Sigma}$ in our error estimates \cref{eq:est_conv_criterion}.} 
    \item \textbf{LOBPCG.}
    Algebraically, apply LOBPCG to the Gram matrix $\mat{A}^* \mat{A}$ with preconditioner $\mat{P}\mat{P}^*$; the implementation never forms $\mat A^* \mat A$ (\cref{sec:lobpcgimplementation}).
    For $i=0,1,2,\ldots$, do the following:
    \begin{enumerate}
        \item \textbf{Generate new search direction.} Set
        \begin{equation} \label{eq:search-direction}
            \vec{w}_i \coloneqq \mat{P}\big(\mat{P}^*\big(\mat{A}^*(\mat{A}\vec{v}_i) - \norm{\mat{A}\vec{v}_i}^2 \cdot \vec{v}_i \big) \big).
        \end{equation}
        \item \textbf{Compute minimum-energy direction in trial space.} Set
        \begin{equation} \label{eq:rr-intro}
            \vec{v}_{i+1} \coloneqq \argmin \{ \norm{\mat{A}\vec{z}} : \vec{z} \in \spn \{ \vec{v}_{i-1},\vec{v}_i,\vec{w}_i \}, \norm{\vec{z}} = 1\} .
        \end{equation}
        This minimization problem can be solved in $\order(mn)$ operations via an SVD of $\mat{A}[\vec{v}_{i-1}, \vec{v}_i, \vec{w}_i]$.
    \end{enumerate}
\end{enumerate}
Under an appropriate condition on the singular value gap and the quality of the random embedding $\mat{S}$, we prove geometric convergence of $\vec{v}_i$ to the right singular vector $\vec{v}_{\rm min}$ of $\mat{A}$ associated with the minimum singular value; empirically, we find it converges under much wider conditions than predicted by theory. For difficult problems with small singular value gaps or where one needs multiple singular vectors, one can use a block iterate $\mat{V}_i \in \complex^{n\times b}$; see \cref{sec:block} for a discussion. 

\subsection{Illustrative example}
To evaluate the RLOBPCG method, we consider a large-scale example that demonstrates its convergence behavior and efficiency relative to other methods without randomized preconditioning. We compare the following methods:
\begin{enumerate}
    \item RLOBPCG with the recommended implementation (\cref{alg:rlopcg}),
    \item LOBPCG with the diagonal preconditioner $\mat{P} = \diag(\norm{\mat{A}(:,i)}^{-1})$,
    \item LOBPCG with no preconditioner,
    \item The Lanczos algorithm \cite[sec.~10.1]{GV13} applied to $\mat{A}^*\mat{A}$ (equivalently, Golub--Kahan bidiagonalization \cite[sec.~10.4]{GV13}),
    \item The sketch-and-solve approximate minimum singular vector $\vec{v}_0$ \cite{PN23}, and
    \item MATLAB's \texttt{svd} routine.
\end{enumerate}
To ensure that all methods are compared under similar conditions, we initialize all of the iterative methods using the sketch-and-solve initialization $\vec{v}_0$. The experiments are intended to illustrate the effects of different preconditioners for LOBPCG, with unpreconditioned Lanczos, sketch-and-solve, and the SVD provided as baselines.

We test two examples, of dimension $m = 2\cdot 10^6$ by $n = 4\cdot 10^3$:
\begin{itemize}
    \item \textbf{Easy singular values.} To create an easy problem $\mat{A} = \mat{U}\diag(\sigma_1,\ldots,\sigma_n)\mat{V}^*$, we generate $\mat{V}$ using the MATLAB command
    \begin{equation*}
        \mat{V} = \texttt{orth(diag(randn(n,1))+1e-4*randn(n))} \in \real^{n\times n},    \end{equation*}
    we let $\mat{U} \in \real^{m\times n}$ be a Haar-random matrix with orthonormal columns, and we set the first $n-1$ singular values $\sigma_1,\ldots,\sigma_{n-1}$ to be geometrically spaced between $1$ and $0.1$ and the smallest singular value to be $\sigma_n \coloneqq 10^{-7}$.
    This problem is relatively easy because there is a large gap between the smallest and second-smallest singular values.
    Also, the right singular vectors are coherent (well-aligned with the standard basis vectors), making diagonal preconditioning fairly effective.
    \item \textbf{Hard singular values.} To generate a hard problem $\mat{A} = \mat{U}\diag(\sigma_1,\ldots,\sigma_n)\mat{V}^*$, we choose both $\mat{U} \in \real^{m\times n}$ and $\mat{V} \in \real^{n\times n}$ to be Haar-random, and we choose the singular values $\sigma_1,\ldots,\sigma_n$ to be geometrically spaced between $1$ and $10^{-10}$.
    This problem is challenging because there is a small relative gap $(\sigma_{n-1}^2 - \sigma_n^2) / \sigma_n^2 \approx 0.01$ at the minimum singular value.
    It is also more difficult for the diagonally preconditioned LOBPCG method, because the problem is incoherent.
\end{itemize}

\begin{figure}[t]
    \centering
    \begin{subfigure}[t]{0.48\textwidth}
        \centering
        \includegraphics[width=\linewidth]{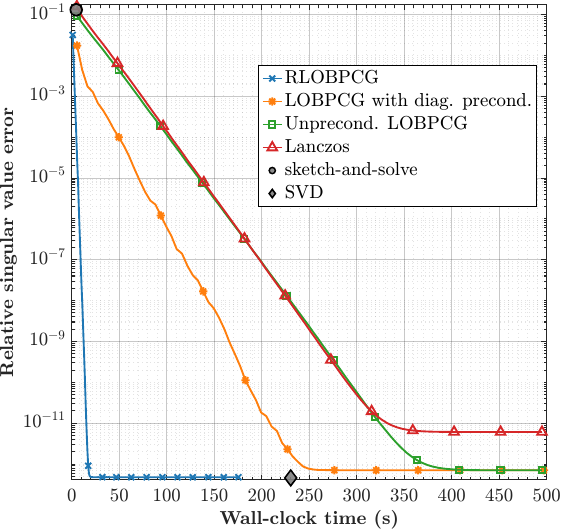}
        \label{subfig:coherenteasyTiming}
    \end{subfigure}%
    \hfill
    \begin{subfigure}[t]{0.48\textwidth}
        \centering
        \includegraphics[width=\linewidth]{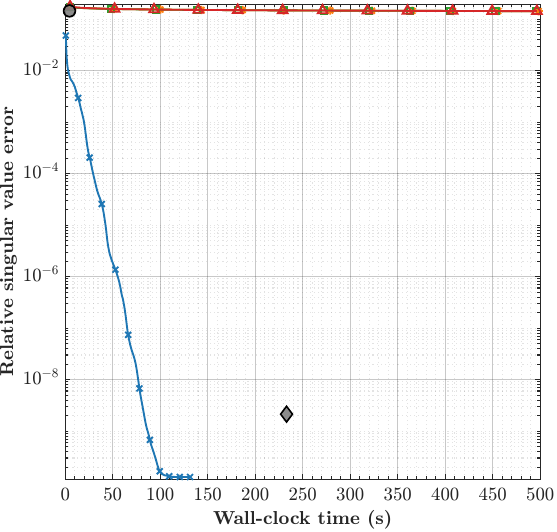}
                \label{subfig:incoherentdifficultTiming}
    \end{subfigure}
    \caption{Relative singular value error $(\norm{\mat{A}\vec{v}_k} - \sigma_n) / \sigma_n$ vs.\ wall-clock time for iterative minimum singular value algorithms on \textbf{easy} (\emph{left}) and \textbf{hard} (\emph{right}) problems.}
    \label{fig:timing}
\end{figure}

\Cref{fig:timing} charts the relative error  in the minimum singular value $(\norm{\mat{A}\vec{v}_k} - \sigma_n) / \sigma_n$ against wall-clock time for each method on both problems.
On the \textbf{easy} problem (\emph{left}), all of the iterative algorithms reach nearly the same accuracy as MATLAB's \texttt{svd} with differing convergence rates.
Even on this easy problem, RLOBPCG is substantially faster than the methods lacking randomized preconditioning, achieving maximal accuracy $\mathbf{12.1\times}$ \textbf{faster} than \texttt{svd}. 
On the \textbf{hard} problem (\emph{right}), all of the non-RLOBPCG solvers fail to make any progress in the allotted time window.
In contrast, RLOBPCG remains robust and rapidly convergent, attaining the same accuracy as the SVD with a $\mathbf{2.8\times}$ \textbf{speedup}.
This example highlights the effectiveness of randomized preconditioning, which preserves convergence even when other methods fail.

\subsection{Summary of theoretical results}
The experiments in the previous section suggest that the RLOBPCG method converges geometrically when there is a singular value gap.
Our theoretical results confirm this to be the case.

To set the stage for our results, we employ the following standard definition:

\begin{definition}[Subspace embedding] \label{def:subspace_embedding}
    A (random) matrix $\mat{S} \in \complex^{d\times m}$ is called a \emph{subspace embedding} for a matrix $\mat{A}\in \complex^{m\times n}$ with \emph{distortion} $0 < \eta < 1$ if
    \begin{equation} \label{eq:distortion_factor}
        (1-\eta)\norm{\vec{Ay}}^2 \leq \norm{\mat{S}\vec{Ay}}^2 \leq (1+\eta) \norm{\vec{Ay}}^2 \quad \text{for all $\vec{y} \in \complex^n$}.
    \end{equation}
\end{definition}

The concept of a subspace embedding is due to Sarl\'os \cite{Sar06}; see \cite[sec.~5]{KT24} for an introduction to idea.
Most constructions of subspace embeddings are randomized and satisfy the condition \cref{eq:distortion_factor} with high probability.
Definitions of subspace embeddings in the literature differ on whether the norms in \cref{eq:distortion_factor} are squared; we adopt the present version as it is most convenient for our setting.

With this context, our main result for RLOBPCG is as follows:

\begin{theorem}[RLOBPCG: Convergence] \label{thm:rlobpcg_convergence}
    Let $\mat{A} \in \complex^{m\times n}$ be a tall matrix with singular values $\sigma_1,\ldots,\sigma_n > 0$, and let $\mat{S} \in \complex^{d\times m}$ be a subspace embedding for $\mat{A}$ with distortion $\eta \in [0,1)$.
    Introduce the relative gap between the \emph{squared} singular values
    \begin{equation*}
        \gap \coloneqq \frac{\sigma_{n-1}^2 - \sigma_n^2}{\sigma_n^2}.
    \end{equation*}
    Suppose that the distortion satisfies 
    \begin{equation} \label{eq:distortion_gap}
        \eta < \frac{\gap}{2+\gap}.
    \end{equation}
    Then 
    $\vec{v}_k$, the approximation to 
    $\vec{v}_{\min}$ (the right singular vector of $\mat{A}$ corresponding to the smallest singular value)
    after $k$ iterations of RLOBPCG, satisfies
    \begin{align}
        \frac{\norm{\mat{A}\vec{v}_k}^2 - \sigma_n^2}{\sigma_n^2} &\le C \cdot \gap \cdot q^{2k}; \label{eq:thm_conclusion_1}\\
        |{\sin \measuredangle (\vec{v}_k,\vec{v}_{\rm min})}| \le |{\tan \measuredangle (\vec{v}_k,\vec{v}_{\rm min})}| &\le  C^{1/2} \cdot q^k,\label{eq:thm_conclusion_2}
    \end{align}
    where the prefactor is $C = 2\eta/[(1-\eta)\gap - 2\eta]$ and the convergence rate is 
    \begin{equation} \label{eq:thm_convergence_rate}
        q = \eta + \frac{1-\eta}{1+\gap}.
    \end{equation}
\end{theorem}

As a corollary, we get a runtime guarantee for RLOBPCG.

\begin{corollary}[RLOBPCG: Runtime] \label{cor:rlobpcg_runtime}
    Set $\eta \coloneqq \min(\gap/3,1/6)$ and let $\varepsilon \in (0,24]$.
    Execute RLOBPCG using a sparse sign embedding $\mat{S} \in \complex^{d\times m}$ with embedding dimension $d = \order((n \log n) / \eta^2)$ and sparsity parameter $\order((\log n)/\eta)$.
    Then, with 99\% probability, RLOBPCG produces a unit vector $\vec{v}_k$ satisfying
    \begin{equation*}
        \norm{\mat{A}\vec{v}_k}^2 \le (1+\varepsilon) \sigma_n^2
    \end{equation*}
    in $\order (\log(1/\varepsilon)/\min(\gap,1))$ iterations    and
    \begin{equation*}
         \order\left(\frac{\operatorname{nnz}(\mat{A})}{\min(\gap,1)} \log \left( \frac{n}{\varepsilon} \right) + \frac{n^3 \log n}{\min(\gap,1)^2}\right) \text{ total operations}.
    \end{equation*}
    Here, $\operatorname{nnz}(\mat{A})$ is the number of nonzero entries in $\mat{A}$.
\end{corollary}
We note several positive features of our result.
First, LOBPCG is typically analyzed from an arbitrary starting vector $\vec{v}_0$, which must satisfy strict requirements to ensure convergence to global minimum.
The RLOBPCG method uses the ``sketch-and-solve'' starting vector $\vec{v}_0$, which satisfies the appropriate conditions with high probability.
Second, when implemented with an embedding of sufficiently small $\eta$ relative to the singular value gap, RLOBPCG produces a high-quality preconditioner which ensures geometric convergence to the minimum singular value at a constant rate.

Let us also draw attention to a basic weakness in our result. To obtain any guarantee, we require the stringent condition \cref{eq:distortion_gap}, which ensures the distortion of the embedding is small relative to the singular value gap. The condition is severe when the gap is extremely small.
For example, when the gap is approximately $10^{-5}$, then \cref{eq:distortion_gap} implies that we need an embedding of distortion $\eta \approx 10^{-5}$, which requires a gigantic embedding dimension of size $d \approx n/\eta^2 \approx 10^{10}n$.
In practice, we have always observed convergence of $\norm{\mat{A}\vec{v}_k}$ (perhaps slowly) to the minimum singular value under the much weaker condition that $\eta$ is a small constant bounded away from $1$ (e.g., $\eta \le 1/2$).
Removing this limitation in our analysis seems challenging; proving global convergence of LOBPCG or any other preconditioned eigensolver from a ``generic'' initialization remains an open problem \cite{KN03,Ovt06,AKNOEZ17,AKSV24}. Improving our analysis would require overcoming this fundamental limitation, showing that LOBPCG (or another eigensolver) converges to the minimum singular value even when $\norm{\mat{A}\vec{v}_0} > \sigma_{n-1}$.

Our analysis of RLOBPCG uses results for preconditioned inverse iteration (PINVIT) \cite{KN03} or preconditioned steepest descent \cite{Ney12}, both of which are guaranteed to be outperformed by LOBPCG for a single step. Our proof uses only this single-step comparison, so \cref{thm:rlobpcg_convergence} is effectively a convergence result for PINVIT or preconditioned steepest descent, with the bound transferred to always-no-worse LOBPCG method.
More discussions and proofs appear in \cref{sec:theory}. 

\subsection{Roadmap}

The rest of the paper is structured as follows.
\Cref{sec:algorithms} discusses implementation of RLOBPCG, choice of embedding, and a block version.
\Cref{sec:theory} discusses theoretical results, and \cref{sec:experiments} provides further experiments.
We conclude in \cref{sec:rational} with an application to computing best rational approximations.

\section{Algorithms and implementation} \label{sec:algorithms}
In this section, we discuss implementation details for RLOBPCG. In \cref{sec:lobpcgimplementation}, we describe the RLOBPCG algorithm and its implementation. \Cref{sec:embedding_choice} provides recommendations for the choice of the subspace embedding.
Finally, in \cref{sec:block}, we discuss block implementations of RLOBPCG.

\subsection{RLOBPCG implementation} \label{sec:lobpcgimplementation}
Pseudocode for RLOBPCG is presented below in \cref{alg:rlopcg}, where MATLAB notation is used for matrix indexing.
As outlined in the introduction, the RLOBPCG method consists of four main steps: (1) generate an embedding, (2) sketch and factorize, (3) define the preconditioner and initialization, and (4) apply LOBPCG to $\mat{A}^*\mat{A}$. Choice of embedding for step (1) is discussed in \cref{sec:embedding_choice}, and steps (2) and (3) are straightforward: We compute the SVD of the sketch $\mat{SA}=\mat{\tilde U\tilde S\tilde V}^*$ and set the preconditioner $\mat{P}=\mat{\tilde V\tilde S}^{-1}$. The initial solution $\vec{v}_0$ is then set to be the smallest right singular vector of the sketch $\mat{SA}$. For the rest of this section, we discuss step (4).

\begin{algorithm}[t]
	\caption{RLOBPCG (single vector)} \label{alg:rlopcg}
	\begin{algorithmic}[1]
		\Require Matrix $\mat{A}\in\real^{m\times n}$, tolerance $\tau \in [0,1)$
		\Ensure Minimum right singular vector $\mat{v} \in \real^{n}$ and singular value $\sigma \in \real_{\geq0}$        \State $\mat{S} \gets \text{$4n\times m$ SparseStack with sparsity $\zeta = 4$}$
        \State $[\sim, \mat{\tilde{\Sigma}}, \mat{\tilde{V}}] \gets \Call{SVD}{\mat{S}\mat{A}}$ \Comment{Sketch and take SVD} \label{line:sketch-svd}
        \State $\mat{P} \gets \smash{\mat{\tilde{V}}\mat{\tilde{\Sigma}}}^{-1}$ \Comment{Build preconditioner}
        \State $\vec{v}_0 \gets \mat{\tilde{V}}(\,:\,,n)$, $\mat{A}_{\vec{v}_0} \gets \mat{A}\vec{v}_0$, $\theta_0 \gets \infty$ \Comment{Sketch-and-solve initialization}
        \State $\mat{r}_0 \gets \mat{A}^* \mat{A}_{\vec{v}_0} - \norm{\mat{A}_{\vec{v}_0}}^2 \vec{v}_0$ \Comment{Initial residual} \label{line:lobpcg-start}
        \State $\vec{x}_0 \gets [\;], \mat{A}_{\vec{x}_0} \gets [\;]$ \Comment{Initialize to empty; cf.\ \cref{eq:x-defining}} \label{line:momentum-initialize}
        \For{$i = 0,1,...,\texttt{max\_iter}-1$} \Comment{Cap total iterations (e.g., $\texttt{max\_iter}=200$)}
            \State $\vec{w}_i \gets \mat{P}\big(\smash{\mat{P}^*} \vec{r}_i\big)$ \Comment{Compute search direction}
            \State $\vec{w}_i \gets \vec{w}_i - [\vec{v}_i,\vec{x}_i]\left([\vec{v}_i,\vec{x}_i]^* \vec{w}_i\right)$ \Comment{Orthogonalize $\vec{w}_i$ against $[\vec{v}_i,\vec{x}_i]$} \label{line:cgs1}
            \State $\vec{w}_i \gets $ \Call{Normalize}{$\vec{w}_i - [\vec{v}_i,\vec{x}_i]([\vec{v}_i,\vec{x}_i]^* \vec{w}_i)$} \Comment{CGS2 for stability} \label{line:cgs2}
            \State $\mat{T} \gets [\vec{v}_i,\vec{x}_i,\vec{w}_i]$ \Comment{Trial space}
            \State $\mat{A}_{\mat{T}} \gets [\mat{A}_{\vec{v}_i},\mat{A}_{\vec{x}_i},\mat{A}\vec{w}_i]$ \Comment{Matrix--vector product with trial space} \label{line:matvecs}
            \State $[\sim,\mat{\Theta},\mat{C}] \gets \Call{SVD}{\mat{A}_{\mat{T}}}$ \Comment{Energy minimization \cref{eq:rr-intro}} \label{line:rayleigh-ritz}
            \State $\vec{v}_{i+1} \gets \mat{T}\mat{C}(\,:\,,3), \theta_{i+1} \gets \mat{\Theta}(3,3)$ \Comment{Updated next iterate}
            \State $\mat{A}_{\vec{v}_{i+1}} \gets \mat{A}_{\mat{T}}\mat{C}(\,:\,,3)$ \Comment{Updated $\mat{A}\vec{v}_{i+1}$}
            \State $\vec{r}_{i+1} \gets \mat{A}^* \mat{A}_{\vec{v}_{i+1}} - \theta_{i+1}^2 \vec{v}_{i+1}$ \Comment{Residual} \label{line:residual}
            \If{$\Call{Mod}{i,5}=0$} \Comment{Check for convergence every five iterations} 
                \State $\mathtt{tol\_check} \gets (\norm{\vec{r}_{i+1}} \leq \mat{\tilde{\Sigma}}(1,1)^2 \cdot\tau)$ \Comment{Tolerance \cref{eq:est_conv_criterion}}
                \State $\mathtt{stag\_check\_1} \gets (1.1 \norm{\vec{r}_{i+1}} \ge \norm{\vec{r}_{i-4}})$
                \State $\mathtt{stag\_check\_2} \gets (1.1 \cdot (\theta_{i-4} - \theta_{i+1})/\theta_i \ge (\theta_{i-9} - \theta_{i-4})/\theta_{i-4})$
                \State \textbf{if} $\mathtt{tol\_check}$ \textbf{and} $\mathtt{stag\_check\_1}$ \textbf{and} $\mathtt{stag\_check\_2}$ \textbf{break} 
            \EndIf
            \State $\vec{q} \gets \Call{Normalize}{\mat{C}(1,1\!:\!2)^*}$
            \State $\mat{x}_{i+1} \gets \mat{T}\mat{C}(:,1\!:\!2)\vec{q}$, $\mat{A}_{\mat{x}_{i+1}} \gets \mat{A}_{\mat{T}}\mat{C}(:,1\!:\!2)\vec{q}$ \Comment{Update $\vec{x}_{i+1},\vec{Ax}_{i+1}$; cf.\ \cref{eq:x-defining}}\label{line:lobpcg-stop}
        \EndFor
        \State $\vec{v} \gets \vec{v}_{i+1}$, $\sigma \gets \theta_{i+1}$ \Comment{Collect outputs}
	\end{algorithmic}
\end{algorithm}
\paragraph{Orthogonalization and minimizing matrix--vector products}
The dominant computational cost in LOBPCG consists of computing matrix--vector products with $\mat{A}$ and orthogonalizing the basis for the energy-minimization step \cref{eq:rr-intro}.
Recall that each LOBPCG iteration computes the minimum-energy direction in a trial space
\begin{equation*}
    \vec{v}_{i+1} = \argmin \{ \norm{\mat{A}\vec{z}} : \vec{z} \in \spn \{ \vec{v}_{i-1},\vec{v}_i,\vec{w}_i \}, \norm{\vec{z}} = 1\}.
\end{equation*}
The search direction $\vec{w}_i$ was defined above in \cref{eq:search-direction}.
A direct implementation would orthogonalize the trial space $[\vec{v}_{i-1},\vec{v}_i,\vec{w}_i]$ and perform Rayleigh--Ritz projection onto this basis at each iteration. This approach would require three matrix--vector products with $\mat{A}$ and an orthogonalization of a three-dimensional basis per iteration.
We can reduce this cost by carefully maintaining and updating the basis, reducing the cost to just one matrix--vector product per iteration \cite{HL06,DSYG18}.

This efficient LOBPCG implementation for $\mat{A}^* \mat{A}$ is given in lines \ref{line:lobpcg-start}--\ref{line:lobpcg-stop} of \cref{alg:rlopcg}.
As the core ``trick'', we introduce a vector $\vec{x}_i$ (line \ref{line:momentum-initialize}) such that
\begin{equation} \label{eq:x-defining}
    \{\vec{v}_i,\vec{x}_i\} \text{ is an orthonormal basis for } \spn \{\vec{v}_{i-1},\vec{v}_i\}.
\end{equation}
 
 Since $\{\vec{v}_i,\vec{x}_i\}$ are orthonormal, only $\vec{w}_i$ needs to be orthogonalized against these two vectors, and only a single matrix--vector product $\mat{A}\vec{w}_i$ is required per iteration (line \ref{line:matvecs}).
 To combine speed and numerical stability, we perform orthogonalization by  classical Gram--Schmidt with reorthogonalization (CGS2, \ref{line:cgs1}--\ref{line:cgs2}).
 Faster but less stable approaches, such as SVQB \cite{SW02,HL06}, are also possible for this step.
 It is important to note that as $\vec{v}_i$ converges to a singular vector, it becomes increasingly linearly dependent with $\vec{v}_{i-1}$. Therefore, the orthogonalization step is crucial---skipping it or using unstable orthogonalization methods may introduce spurious singular values \cite{Kny01}.
 The energy-minimization step \cref{eq:rr-intro} is performed by an SVD of $\mat A\mat T$ (line \ref{line:rayleigh-ritz}), rather than an eigendecomposition of $(\mat A \mat T)^*(\mat A \mat T)$. With an orthonormal trial basis $\mat T$, this extraction gives the Rayleigh--Ritz step for LOBPCG on $\mat A^* \mat A$. In exact arithmetic, the single-vector RLOBPCG iteration is closely related to the Golub--Kahan Davidson method with $+1$ locally optimal restarting and a maximum basis dimension of three~\cite{GSR19}. Our implementation forms the normal equation residual for subspace expansion, so this connection does not imply identical finite--precision behavior.
 A formal stability analysis is a natural subject for future work.
 
\paragraph{Stopping criteria}

Stopping a preconditioner eigensolver for a minimum singular value problem with a small singular value gap is challenging, and we are unaware of a fully satisfactory solution.
For our implementations, we the following ad hoc strategy.
As a first stopping criterion for RLOBPCG, we employ a variant of the following backward error-based convergence criterion from \cite{DSYG18}:
\begin{equation} \label{eq:conv_criterion}
    \frac{\norm{\vec{r}_i}}{\big(\smash{\norm{\mat{A}}}^2+\smash{\norm{\mat{Av}_i}}^2\big)} \leq \tau,
\end{equation} where $\vec{r}_i = \mat{A}^* \mat{A} \vec{v}_i - \norm{\mat{Av}_i}^2 \vec{v}_{i}$ is the residual and $\tau \in [0,1)$ is a tolerance parameter.
To implement \cref{eq:conv_criterion}, we need to compute or estimate $\norm{\mat{A}}$.
By the subspace embedding property (\cref{def:subspace_embedding}), the sketch $\mat{SA}$ provides a good approximation to the singular values of $\mat{A}$, so that $\norm{\mat{A}} \approx \norm{\mat{SA}} = \mat{\tilde{\Sigma}}(1,1)$ for $\mat{\tilde{\Sigma}}$ computed in line \ref{line:sketch-svd} of \cref{alg:rlopcg}.
Using this estimate with the bound $\norm{\mat{Av}_i} \leq \norm{\mat{A}}$, we recommend the following convergence criterion:
\begin{equation} \label{eq:est_conv_criterion}
    \norm{\vec{r}_i} \leq \mat{\tilde{\Sigma}}(1,1)^2\cdot \tau.
\end{equation}
A similar sketch-based estimate of $\norm{\mat{A}}$ is also used in \cite{DSYG18}.

If one wants the maximum-achievable accuracy, we have not found the stopping rule \cref{eq:est_conv_criterion} to be reliable on its own, and we were unable to find a way of setting the tolerance $\tau$ that reliably ensured the method both always halted and achieved near-optimal accuracy.
Thus, we recommend using the stopping rule \cref{eq:est_conv_criterion} with the conservative tolerance $\tau = 2n^{1/2}\varepsilon_{\rm mach}$ and supplementing with two additional checks for stagnation, performed every five iterations:
\begin{equation} \label{eq:stagnation-check}
    1.1 \norm{\vec{r}_i} \ge \norm{\vec{r}_{i-5}} \quad \text{and} \quad
    1.1 \cdot (\theta_{i-5} - \theta_i)/\theta_i \ge (\theta_{i-10} - \theta_{i-5})/\theta_{i-5}.
\end{equation}
Here, $\theta_i$ is the approximation to the minimum singular value produced at iteration $i$.
The method is halted when both the tolerance \cref{eq:est_conv_criterion} and the stagnation checks \cref{eq:stagnation-check} are satisfied.
Even with this two-stage strategy, we have still found that this stopping criteria can be too aggressive, stopping the method when a few more digits of numerical accuracy can still be squeezed out.
The development of more robust stopping criteria is a natural subject for future work.

\subsection{Choice of embedding} \label{sec:embedding_choice}

The embedding needs to provide a preconditioner spectrally equivalent to $(\mat A^*\mat A)^{-1}$. By \cref{prop:randomized_preconditioning} below, a subspace embedding with distortion $\eta$ ensures $\kappa_2(\mat P^* \mat A^* \mat A \mat P) \leq (1+\eta)/(1-\eta)$, showing that the preconditioner quality is \emph{independent} of both the condition number $\kappa_2(\mat A)$ of the unpreconditioned matrix and the spectral gap. This is precisely the type of spectral equivalence required by the classical theory of preconditioned eigensolvers \cite{Ney01,Ney12,KN03}.

Accordingly, RLOBPCG can be implemented with any matrix $\mat{S} \in \complex^{d\times m}$ satisfying the subspace embedding property (\cref{def:subspace_embedding}) with distortion $\eta$ bounded away from $1$; see \cite[secs.~8--9]{MT20} and \cite[chs.~2, 6, \& 7]{MDM+23} for a survey of constructions.
Based on the testing in \cite{DM23,Epp25a,CEMT25}, we recommend using SparseStack embeddings:
\begin{equation*}
    \mat{S} = \zeta^{-1/2} \begin{bmatrix}
        \vec{s}_1 & \cdots & \vec{s}_m
    \end{bmatrix},
\end{equation*}
where each column $\vec{s}_i$ is divided into $\zeta$ equally sized pieces, each of which has a single nonzero entry with uniformly random value $\pm 1$ placed in uniformly random position.
The parameter $1\le \zeta \le d$ controls the level of sparsity.

How should one pick the parameters $d$ and $\zeta$?
Seminal analysis by Cohen \cite{Coh16} yields the following choices:
\begin{fact}[SparseStack] \label{fact:sketching_constructions}
    Fix a matrix $\mat{A} \in \real^{m\times n}$ and $\eta \in (0,1)$.
    A SparseStack is a subspace embedding for $\mat{A}$ for some $d = \order(n\log(n) / \eta^2)$ and sparsity $\zeta = \order(\log(n)/\eta)$ with distortion $\eta$ with 99\% probability.
\end{fact}
Analyses of different sparse dimensionality reduction maps with a different dependency of $d$ and $\zeta$ on $n$ and $\eta$ appear in \cite{CDD24a,CDD25b,Tro25,CEMT25}.
In practice, we use parameters
\begin{equation*}
    2n \le d \le 4n \quad \text{and} \quad \zeta = 4,
\end{equation*}
which we find reliably yields an embedding of sufficient quality for our purposes; \cref{subsubsec:sketchsize} justifies this empirically.

\subsection{Block version} \label{sec:block}
If we are interested in obtaining a minimum right singular subspace rather than a single minimum right singular vector, we turn to block methods.
Block methods are particularly useful when the minimum right singular vector is ill-conditioned because the gap between the smallest and second smallest singular values is very small.
In such cases, projection-based methods struggle to isolate the target (minimum) singular vector \cite{Nak17,Nak20ritz}, and block methods targeting a singular subspace are more appropriate.
Block methods are also more robust to singular value clustering \cite{Ovt06,ZKN25} and make use of hardware-efficient BLAS3 matrix operations.
\begin{algorithm}[t]
	\caption{RLOBPCG (block version)} \label{alg:rlobpcg}
	\begin{algorithmic}[1]
		\Require Matrix $\mat{A}\in\real^{m\times n}$, block size $b$, tolerance $\tau \in [0,1)$
		\Ensure Minimum right singular subspace $\mat{V} \in \real^{n\times b}$ and diagonal matrix containing the minimum singular values $\mat{\Sigma} \in \real^{b\times b}$
        \State $\mat{S} \gets \text{$4n\times m$  SparseStack with sparsity $\zeta = 4$}$
        \State $[\sim, \mat{\tilde{\Sigma}}, \mat{\tilde{V}}] \gets \Call{SVD}{\mat{S}\mat{A}}$ \Comment{Sketch and take SVD}
        \State $\mat{P} \gets \smash{\mat{\tilde{V}}\mat{\tilde{\Sigma}}}^{-1}$ \Comment{Build preconditioner}
        \State $\mat{V}_0 \gets \mat{\tilde{V}}(\;:\;,n\!-\!b\!+\!1\!:\!n)$, $\mat{A}_{\mat{V}_0} \gets \mat{A}\mat{V}_0$ \Comment{Sketch-and-solve initialization} \label{line:block-initial}
        \State $[\mat{C},\mat{\Theta}] \gets \Call{RayleighRitz}{\mat{A}_{\vec{V}_0},\vec{V}_0}$
        \State $\mat{V}_0 \gets \mat{V}_0 \mat{C}$, $\mat{A}_{\vec{V}_0} \gets \mat{A}_{\vec{V}_0} \mat{C}$
        \State $\mat{R}_0 \gets \mat{A}^* \mat{A}_{\vec{V}_0} -  \vec{V}_0 \mat{\Theta}^2$ \label{line:block-initial-residual}
        \State $\vec{X}_0 \gets [\;], \mat{A}_{\vec{X}_0}\gets [\;]$ \Comment{Initialize to empty}
        \For{$i = 0,1,...,\texttt{max\_iter}-1$} \Comment{Cap total iterations (e.g., $\texttt{max\_iter}=100$)}
            \State $\vec{W}_i \gets \mat{P}\big(\smash{\mat{P}^* \vec{R}_i}\big)$ 
            \State $\vec{W}_i \gets \vec{W}_i - [\vec{V}_i,\vec{X}_i]\left([\vec{V}_i,\vec{X}_i]^* \vec{W}_i\right)$ \label{line:block-orth-1}\Comment{Orthogonalize $\mat{W}_i$ against $[\vec{V}_i,\vec{X}_i]$}
            \State $\vec{W}_i \gets$ \Call{Orth}{$\vec{W}_i - [\vec{V}_i,\vec{X}_i]([\vec{V}_i,\vec{X}_i]^* \vec{W}_i)$} \Comment{BCGS2 for stability} \label{line:block-orth-2}
            \State $\mat{T} \gets [\vec{V}_i,\vec{X}_i,\vec{W}_i]$ \Comment{Trial space}
            \State $\mat{A}_{\mat{T}} \gets [\mat{A}_{\vec{V}_i},\mat{A}_{\vec{X}_i},\mat{A}\vec{W}_i]$ \Comment{Matrix--vector products with trial space}
            \State $[\sim,\mat{\Theta},\mat{C}] \gets \Call{SVD}{\mat{A}_{\mat{T}}}$ \Comment{Rayleigh--Ritz}
            \State $\vec{V}_{i+1} \gets \mat{T}\mat{C}(\,:\,,2b\!+\!1\!:\!3b)$ \Comment{Updated next iterate}
            \State $\mat{A}_{\vec{V}_{i+1}} \gets \mat{A}_{\mat{T}}\mat{C}(\,:\,,2b\!+\!1\!:\!3b)$ \Comment{Updated $\mat{A}\vec{V}_{i+1}$}
            \State $\vec{R}_{i+1} \gets \mat{A}^* \mat{A}_{\vec{V}_{i+1}} - \vec{V}_{i+1} \mat{\Theta}(2b\!+\!1\!:\!3b,2b\!+\!1\!:\!3b)^2 $ \Comment{Residual matrix}\label{line:block-residual}
            \State \textbf{if} $\max\limits_{1\leq j \leq b}\norm{\mat{R}_{i+1}(\,:\,,j)} \leq \mat{\tilde{\Sigma}}(1,1)^2 \cdot \tau$ \textbf{then break} \Comment{Better: also use \cref{eq:stagnation-check}}    \label{line:block-stopping}   
            \State $\mat{Q} \gets \Call{Orth}{\mat{C}(1\!:\!b,1\!:\!2b)^*}$
            \State $\mat{X}_{i+1} \gets \mat{T}\mat{C}(\,:\,,1\!:\!2b)\mat{Q}$, $\mat{A}_{\mat{X}_{i+1}} \gets \mat{A}_{\mat{T}}\mat{C}(\,:\,,1\!:\!2b)\mat{Q}$        \EndFor
        \State $\vec{V} \gets \vec{V}_{i+1}$
        \State $\mat{\Sigma} \gets \mat{\Theta}(2b\!+\!1\!:\!3b,2b\!+\!1\!:\!3b)$
	\end{algorithmic}
\end{algorithm}

A block version of RLOBPCG is presented in \cref{alg:rlobpcg}. It is structurally similar to the single vector version (\cref{alg:rlopcg}), with the main change being the move from $1$-dimensional to a $b$-dimensional subspace for a user-specified block size $b$. Similarly, in exact arithmetic, the block iteration is closely related to a block extension of Golub--Kahan Davidson with $+b$ locally optimal restarting and a maximum trial-subspace dimension of $3b$~\cite{GSR19}.

The working block size $b$ need not equal the number of requested singular triplets. \cref{alg:rlobpcg} is written for computing $b$ singular triplets. In the experiments of \cref{sec:blocktest}, only the smallest singular triplet is requested. All $b$ Ritz vectors are retained in the iteration, but the convergence criteria \eqref{eq:est_conv_criterion}--\eqref{eq:stagnation-check} are applied only to the smallest Ritz pair.
Instead of a vector initialization, we now provide an initial subspace using the sketch in line \ref{line:block-initial}, which serves as a crude-but-effective approximation to the minimum $b$ right singular vectors of $\mat{A}$ \cite{PN23}.
In lines \ref{line:block-initial-residual} and \ref{line:block-residual}, the residual is now a matrix whose columns are the residual vectors for each of the $b$ minimum singular vectors we wish to approximate.
In lines \ref{line:block-orth-1}--\ref{line:block-orth-2}, we apply \emph{block} CGS2 (BCGS2) to orthogonalize the matrix search direction $\mat{W}_i$ against both $\mat{V}_i$ and $\mat{X}_i$.
We stop the algorithm using a block modification of the stopping rule \cref{eq:est_conv_criterion}; see line~\ref{line:block-stopping}.
For the best reliability, one can also use the stagnation checks \cref{eq:stagnation-check}.
A more efficient block version would deflate the converged singular vectors, a process known as locking \cite{Kny01,KALO07}. Once an index $j$ meets the condition \cref{eq:est_conv_criterion}, it would be locked by storing it in $\mat{V}_i$, but removing the corresponding columns from $\mat{X}_i$ and $\mat{W}_i$. This ensures the trial space is orthogonal to converged singular vectors. Locking is not used in our implementation. For details on locking, see~\cite[sec.~2.3]{KALO07}.

\section{Theory} \label{sec:theory}

In this section, we present convergence analysis of RLOBPCG.
\Cref{sec:sketching-background,sec:background-lobpcg} provide background on the analysis of sketching and LOBPCG, \cref{sec:interpreting-lobpcg-bounds} discusses how to relate quantities in the LOBPCG convergence theory to the accuracy of singular values and singular vectors, and \cref{sec:rlobpcg_convergence-proof,sec:rlobpcg_runtime-proof} prove \cref{thm:rlobpcg_convergence,cor:rlobpcg_runtime}.

\subsection{Background: Results for sketching} \label{sec:sketching-background}

The definition of subspace embeddings appeared above in \cref{def:subspace_embedding}.
This section will review two necessary results about sketching: the quality of the randomized preconditioner $\mat{P}$ and the quality of the ``sketch-and-solve'' initialization $\vec{v}_0$. 
For the former, we employ the following result of Rokhlin \& Tygert \cite{RT08}; see \cite[Prop.~5.4]{KT24} for a modern proof.

\begin{fact}[Randomized preconditioning] \label{prop:randomized_preconditioning}
    Let $\mat{A} \in \complex^{m\times n}$ be a full-rank matrix, and let $\mat{S} \in \complex^{d\times m}$ be a subspace embedding for $\mat{A}$ of distortion $\eta \in [0,1)$.
    Introduce the preconditioner $\mat{P} = \mat{V} \mat{\Sigma}^{-1}$ from an SVD $\mat{S}\mat{A} = \mat{U}\mat{\Sigma}\mat{V}^*$.
    Then 
    \begin{equation*}
        \frac{1}{1+\eta} \le \sigma_{\min}(\mat{A}\mat{P})^2 \le \sigma_{\rm max}(\mat{A}\mat{P})^2 \le \frac{1}{1-\eta}.
    \end{equation*}
\end{fact}
To quantify the quality of the ``sketch-and-solve'' initialization $\vec{v}_0 = \mat{P}\evec_n/\norm{\mat{P}\evec_n}$ for \sloppy RLOBPCG, we use the following result:
\begin{fact}[Sketch-and-solve for minimum singular value] \label{prop:sketch_and_solve}
    Import the setting of \cref{prop:randomized_preconditioning} and define $\vec{v}_0 \coloneqq \mat{P}\evec_n /\norm{\mat{P}\evec_n}$.
    Then 
    \begin{equation*}
        \norm{\mat{A}\vec{v}_0}^2 \le \frac{1+\eta}{1-\eta} \cdot \sigma_{\rm min}(\mat{A})^2.
    \end{equation*}
\end{fact}
Versions of this result appear in \cite[Rem.~2.1]{PN23}, \cite[Prop.~5.5]{KT24}; see also \cite{GPW12}.
\subsection{Background: Preconditioned eigensolvers} \label{sec:background-lobpcg}

Despite decades of interest, a sharp convergence theory for LOBPCG remains elusive.
Therefore, to understand LOBPCG, it is common to compare the iterates of LOBPCG to other preconditioned eigensolvers, such as preconditioned inverse iteration (PINVIT, \cite{Ney01}) or preconditioned steepest descent \cite{Ney12}, that are easier to analyze.
We shall make use of this standard strategy to analyze RLOBPCG.

Given the current iterate $\vec{v} = \vec{v}_i$ and the previous iterate $\vec{v}_{\rm old} = \vec{v}_{i-1}$, recall that the next iterate generated by LOBPCG is given by the equations
\begin{align}
    \vec{w} &\coloneqq \mat{P}\big(\mat{P}^*\big(\mat{A}^*(\mat{A}\vec{v}) - \norm{\mat{A}\vec{v}}^2 \cdot \vec{v} \big) \big), \label{eq:search_direction} \\
    \vec{v}_{\rm LOBPCG} &\coloneqq \argmin \{ \norm{\mat{A}\vec{z}} : \vec{z} \in \spn \{ \vec{v}_{\rm old},\vec{v},\vec{w} \}, \norm{\vec{z}} = 1\} .\label{eq:lobpcg_step}
\end{align}
The next iterate $\vec{v}_{\rm LOBPCG}$ generated by LOBPCG depends
on both the previous iterates $\vec{v}$ and $\vec{v}_{\rm old}$ and the search direction $\vec{w}$.
The preconditioned steepest descent iteration is the same as LOBPCG, except the dependence on the past state $\vec{v}_{\rm old}$ is removed:
\begin{equation}
    \vec{v}_{\rm PSD} \coloneqq \argmin \{ \norm{\mat{A}\vec{z}} : \vec{z} \in \spn \{ \vec{v},\vec{w} \}, \norm{\vec{z}} = 1\} . \label{eq:psd_step}
\end{equation}
Since LOBPCG searches over a larger set of possible vectors, the norm $\norm{\mat{A}\vec{v}}$ is always smaller for LOBPCG than for preconditioned steepest descent:
\begin{equation} \label{eq:lobpcg_psd_compare}
    \norm{\mat{A}\vec{v}_{\rm LOBPCG}} \le \norm{\mat{A}\vec{v}_{\rm PSD}}.
\end{equation}
The relation \cref{eq:lobpcg_psd_compare} allows us to transfer single-step bounds for preconditioned steepest descent to single-step bounds for LOBPCG.

We import the following standard bound for PINVIT with an optimal damping factor \cite[Thm.~1.1]{Ney12} and specialize to our setting. 

\begin{fact}[Preconditioned steepest descent]
    Let $\vec{v} \in \complex^n$ be a vector, let $\mat{A} \in \complex^{m\times n}$ be a full-rank matrix with singular values $\sigma_1,\ldots,\sigma_n$, and let $\mat{P} \in \complex^{n\times n}$ be a full-rank preconditioner.
    Introduce the parameter
    \begin{equation} \label{eq:gamma}
        \gamma \coloneqq \inf_{t > 0} \norm{\Id - t\, \mat{P}^*\mat{A}^*\mat{A}\mat{P}}.
    \end{equation}
    Assume that the unit vector $\vec{v}$ satisfies $\sigma_i \le \norm{\mat{A}\vec{v}} < \sigma_{i-1}$ and introduce the gap
    \begin{equation*}
        \gap_i \coloneqq \frac{\sigma_{i-1}^2 - \sigma_i^2}{\sigma_i^2}.
    \end{equation*}
    Then the preconditioned steepest descent iterate \cref{eq:psd_step} either is below the $i$th singular value $\norm{\mat{A}\vec{v}_{\rm PSD}} \le \sigma_i$ or it admits the bound
    \begin{equation*}
        \frac{\norm{\mat{A}\vec{v}_{\rm PSD}}^2 - \sigma_i^2}{\sigma_{i-1}^2 - \norm{\mat{A}\vec{v}_{\rm PSD}}^2} \le q^2 \cdot \frac{\norm{\mat{A}\vec{v}}^2 - \sigma_i^2}{\sigma_{i-1}^2 - \norm{\mat{A}\vec{v}}^2} \quad \text{for } q = \gamma + \frac{1-\gamma}{1+\gap_i}.
    \end{equation*}
\end{fact}

Applying the comparison \cref{eq:lobpcg_psd_compare}, specializing to the minimum singular value (i.e., $i=n$), and iterating, we obtain a bound for LOBPCG:
\begin{equation} \label{eq:lobpcg_iteration_bound}
    \frac{\norm{\mat{A}\vec{v}_k}^2 - \sigma_n^2}{\sigma_{n-1}^2 - \norm{\mat{A}\vec{v}_k}^2} \le q^{2k} \cdot \frac{\norm{\mat{A}\vec{v}_0}^2 - \sigma_n^2}{\sigma_{n-1}^2 - \norm{\mat{A}\vec{v}_0}^2} \quad \text{provided } \norm{\mat{A}\vec{v}_0} < \sigma_{n-1}.
\end{equation}

\subsection{Interpreting \cref{eq:lobpcg_iteration_bound}} \label{sec:interpreting-lobpcg-bounds}

The bound \cref{eq:lobpcg_iteration_bound} may not initially seem very transparent.
Before proving our main theorem, let us simplify.
First, we bound the approximate minimum singular value $\norm{\mat{A}\vec{v}}$, then we treat the angle $\measuredangle(\vec{v},\vec{v}_{\rm min})$ with the minimum singular vector.
These results are more-or-less standard, and we review them for completeness.

We begin with the approximate minimum singular value $\norm{\mat{A}\vec{v}}$.
Introduce the parameter
\begin{equation*}
    \gap \coloneqq \gap_n = \frac{\sigma_{n-1}^2 - \sigma_n^2}{\sigma_n^2}.
\end{equation*}
as in \cref{thm:rlobpcg_convergence}, and write $B$ for the right-hand side of \cref{eq:lobpcg_iteration_bound}.
Using these notations, \cref{eq:lobpcg_iteration_bound} can be written 
\begin{equation*}
    \norm{\mat{A}\vec{v}_k}^2 - \sigma_n^2 \le B\cdot ((1+\gap)\sigma_n^2 - \norm{\mat{A}\vec{v}_k}^2).
\end{equation*}
Rearranging gives
\begin{equation*}
    \frac{\norm{\mat{A}\vec{v}_k}^2 - \sigma_n^2}{\sigma_n^2} \le \frac{B}{1+B} \cdot \gap \le B \cdot \gap.
\end{equation*}
Reinstating the definition of $B$ as the right-hand side of \cref{eq:lobpcg_iteration_bound}, we obtain
\begin{equation} \label{eq:singular_value_bound}
    \frac{\norm{\mat{A}\vec{v}_k}^2 - \sigma_n^2}{\sigma_n^2} \le q^{2k} \cdot \frac{\norm{\mat{A}\vec{v}_0}^2 - \sigma_n^2}{\sigma_{n-1}^2 - \norm{\mat{A}\vec{v}_0}^2} \cdot \gap \quad \text{provided } \norm{\mat{A}\vec{v}_0} < \sigma_{n-1}.
\end{equation}
The bound \cref{eq:singular_value_bound} may be easier to interpret than the original expression \cref{eq:lobpcg_iteration_bound}.

Now we examine the angle $\measuredangle(\vec{v},\vec{v}_{\rm min})$.
Any unit vector $\vec{v}$ we may written as a linear combination
\begin{equation*}
    \vec{v} = \alpha\, \vec{v}_{\rm min} + \beta\, \vec{v}_\perp \quad \text{with } |\alpha|^2 + |\beta|^2 = 1
\end{equation*}
of $\vec{v}_{\rm min}$ and a unit vector $\vec{v}_\perp$ orthogonal to $\vec{v}_{\rm min}$.
The sine of the angle between $\vec{v}$ and $\vec{v}_{\rm min}$ is $|{\sin \measuredangle(\vec{v},\vec{v}_{\rm min})}| = |\beta|$.
We have the decomposition
\begin{equation*}
    \norm{\mat{A}\vec{v}}^2 = |\alpha|^2 \sigma_n^2 + |\beta|^2 \norm{\mat{A}\vec{v}_\perp}^2 = (1-\sin^2 \measuredangle(\vec{v},\vec{v}_{\rm min}))\sigma_n^2 + \sin^2 \measuredangle(\vec{v},\vec{v}_{\rm min})\norm{\mat{A}\vec{v}_\perp}^2.
\end{equation*}
Since $\vec{v}_\perp$ is orthogonal to $\vec{v}_{\rm min}$, $\norm{\mat{A}\vec{v}_\perp}^2 \ge \sigma_{n-1}^2$.
Thus,
\begin{equation*}
    \norm{\mat{A}\vec{v}}^2 \ge (1-\sin^2 \measuredangle(\vec{v},\vec{v}_{\rm min}))\sigma_n^2 + \sin^2 \measuredangle(\vec{v},\vec{v}_{\rm min})\sigma_{n-1}^2.
\end{equation*}
Rearranging yields
\begin{equation*}
    \sin^2 \measuredangle(\vec{v},\vec{v}_{\rm min}) \le \frac{\norm{\mat{A}\vec{v}}^2 - \sigma_n^2}{\sigma_{n-1}^2-\sigma_n^2}.
\end{equation*}
Thus
\begin{equation*}
    \tan^2 \measuredangle(\vec{v},\vec{v}_{\rm min}) = \frac{\sin^2 \measuredangle(\vec{v},\vec{v}_{\rm min})}{1-\sin^2 \measuredangle(\vec{v},\vec{v}_{\rm min})} \le \frac{\tfrac{\norm{\mat{A}\vec{v}}^2 - \sigma_n^2}{\sigma_{n-1}^2-\sigma_n^2}}{1-\tfrac{\norm{\mat{A}\vec{v}}^2 - \sigma_n^2}{\sigma_{n-1}^2-\sigma_n^2}} = \frac{\norm{\mat{A}\vec{v}}^2 - \sigma_n^2}{\sigma_{n-1}^2-\norm{\mat{A}\vec{v}}^2}.
\end{equation*}
Combining this bound with \cref{eq:lobpcg_iteration_bound} yields
\begin{equation}\label{eq:singular_vector_bound}
    |{\tan \measuredangle(\vec{v},\vec{v}_{\rm min})}| \le q^k \cdot \left(\frac{\norm{\mat{A}\vec{v}_0}^2 - \sigma_n^2}{\sigma_{n-1}^2 - \norm{\mat{A}\vec{v}_0}^2}\right)^{1/2} \quad \text{provided } \norm{\mat{A}\vec{v}_0} < \sigma_{n-1}.
\end{equation}

\subsection{Proof of \cref{thm:rlobpcg_convergence}} \label{sec:rlobpcg_convergence-proof}

In order to instantiate the bound \cref{eq:lobpcg_iteration_bound}, we must bound the parameter $\gamma$ in \cref{eq:gamma} and ensure that $\norm{\mat{A}\vec{v}_0} < \sigma_{n-1}$.

First, we bound $\gamma$.
By \cref{prop:randomized_preconditioning}, the squared singular values of $\mat{A}\mat{P}$ lie between $1/(1+\eta)$ and $1/(1-\eta)$.
Ergo, 
\begin{equation*}
    \gamma = \inf_{t>0} \norm{\Id - t \mat{P}^*\mat{A}^*\mat{A}\mat{P}} \le \inf_{t>0} \max \left\{ \frac{t}{1-\eta}-1, 1 - \frac{t}{1+\eta} \right\}.
\end{equation*}
The right-hand side is minimized by taking $t = 1-\eta^2$, yielding the bound $\gamma \le \eta$.

Next, we establish $\norm{\mat{A}\vec{v}_0} < \sigma_{n-1}$.
By \cref{prop:sketch_and_solve}, the $\vec{v}_0$ satisfies
\begin{equation} \label{eq:sketch_and_solve}
    \norm{\mat{A}\vec{v}_0}^2 \le \frac{1+\eta}{1-\eta} \cdot \sigma_n^2.
\end{equation}
The right-hand side is strictly less than $\sigma_{n-1}^2$ under the hypothesis \cref{eq:distortion_gap}.

Finally, we bound the quantity 
\begin{equation*}
    \frac{\norm{\mat{A}\vec{v}_0}^2 - \sigma_n^2}{\sigma_{n-1}^2 - \norm{\mat{A}\vec{v}_0}^2}
\end{equation*}
Using the bound \cref{eq:sketch_and_solve} and the definition $\sigma_{n-1}^2 = (1+\gap) \sigma_n^2$, we have
\begin{equation*}
    \frac{\norm{\mat{A}\vec{v}_0}^2 - \sigma_n^2}{\sigma_{n-1}^2 - \norm{\mat{A}\vec{v}_0}^2} \le \frac{\tfrac{1+\eta}{1-\eta} - 1}{1 + \gap - \tfrac{1+\eta}{1-\eta}} = \frac{2\eta}{(1-\eta)\gap - 2\eta}.
\end{equation*}
The conclusions \cref{eq:thm_conclusion_1,eq:thm_conclusion_2} of the theorem follow from this equation and the bounds \cref{eq:singular_value_bound,eq:singular_vector_bound} from the previous section.

\subsection{Proof of \cref{cor:rlobpcg_runtime}}\label{sec:rlobpcg_runtime-proof}
Since $\eta = \min(\gap/3,1/6)$, the convergence rate $q$ from \cref{eq:thm_convergence_rate} satisfies the bound $q \le 1 - \min(\gap,1)/6$.Further, we have the bound $C\cdot \gap \le 2$.The claimed complexity results follow from \cref{thm:rlobpcg_convergence} and \cref{fact:sketching_constructions}.

\section{Experiments}
\label{sec:experiments}
This section presents a set of numerical experiments designed to evaluate the accuracy, robustness, and efficiency of RLOBPCG (\cref{alg:rlopcg}) for computing the smallest singular value and its associated singular vector.
All experiments except the block experiment in \cref{sec:blocktest} were performed in MATLAB version R2024b on a workstation with Intel Xeon Gold 6538Y+ CPU @ 2.40GHz (128 cores) and 2TB memory. The block experiment was performed in MATLAB version R2026a with two AMD EPYC 7742 processors (128 cores) and approximately 540 GB memory.

Each curve will be divided into an \emph{active} and a \emph{inactive} segment.
The active portion (colored) corresponds to iterations before which RLOBPCG the convergence criteria \cref{eq:est_conv_criterion,eq:stagnation-check} have not been met.
Once these criteria have been met, the curve transitions to the gray inactive part.
For a given approximation $\vec{\hat{v}}_{\rm min}$ to the minimum singular vector and its corresponding singular value approximation $\hat{\sigma}_{\rm min}$, we will plot both the singular value error $|\hat{\sigma}_{\rm min} - \sigma_{\rm min}| / \sigma_{\rm min}$ and the angle $\sin \angle(\vec{v}_{\rm min},\vec{\hat{v}}_{\rm min})$ to the true minimum singular vector.
In all plots, dashed horizontal lines indicate the level of accuracy predicted by matrix perturbation theory \cite{DK70,Wed72} for a backward stable algorithm.
For the relative singular value $|\hat{\sigma}_{\rm min} - \sigma_{\rm min}| / \sigma_{\rm min}$, the expected error is $u \cdot \sigma_{\rm max} / \sigma_{\rm min}$ and for the angle $\sin \angle(\vec{v}_{\rm min},\vec{\hat{v}}_{\rm min})$, the expected error is $u/\mathrm{gap}_{\rm abs}$.\footnote{The relative gap in the convergence theorem governs convergence in exact arithmetic whereas $\mathrm{gap}_{\rm abs} = (\sigma_{n-1}-\sigma_n)/\sigma_1$ here describes sensitivity to perturbations. The reference levels in the experiments concern accuracy in the presence of roundoff errors, rather than the exact arithmetic convergence rate. Indeed, if the computed triplet is exact for a backward perturbation $\mat A + \mat E$ with $\|\mat E\| \lesssim u\norm{\mat A}$, then Weyl's theorem implies $|\hat{\sigma}_{\min} - \sigma_{\min}| \lesssim u \sigma_{\max}$ and Wedin's theorem~\cite{Wed72} gives $\sin \angle(\vec{v}_{\rm min},\vec{\hat{v}}_{\rm min}) \lesssim \frac{u}{\mathrm{gap}_{\rm abs}}$. The symbol $\lesssim$ in both bounds suppresses $\poly(n)$ prefactors, so agreement within an order of magnitude is the natural benchmark.}
Here, $u\approx 10^{-16}$ is the unit roundoff in double precision and $\mathrm{gap}_{\rm abs} = (\sigma_{n-1} - \sigma_n)/\sigma_1$ is the singular value gap.

\subsection{Single-vector tests}

In this section, we examine how the spectral gap and the condition number affect the performance of RLOBPCG. Throughout this subsection, 
$\mat{A}$ is a $10^4 \times 10^2$ matrix with left and right singular vectors drawn from the Haar distribution.

\subsubsection{Effect of the singular value gap}

First, we investigate how the size of the (relative) singular value gap affects the convergence behavior of RLOBPCG. The singular values decay geometrically from $1$ to $10^{-9}$ for the first $n-2$ values.
We set the smallest singular value to $\sigma_n = 10^{-10}$, and choose the second-smallest singular value $\sigma_{n-1}$ to achieve a prescribed value of the gap, defined in \cref{sec:theory} as
\begin{equation*}
    \mathrm{gap} = \frac{\sigma_{n-1}^2 - \sigma_n^2}{\sigma_n^2}.
\end{equation*}

\begin{figure}[t]
    \centering
    \begin{subfigure}[t]{0.48\textwidth}
        \centering
        \includegraphics[width=\linewidth]{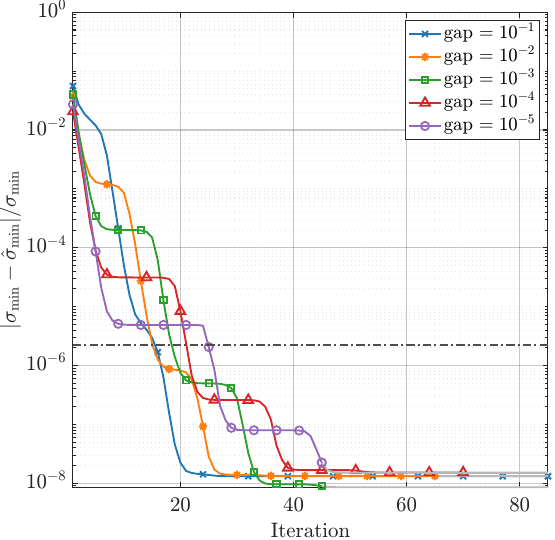}        \label{subfig:gapSingval}
    \end{subfigure}%
    \hfill
    \begin{subfigure}[t]{0.48\textwidth}
        \centering
        \includegraphics[width=\linewidth]{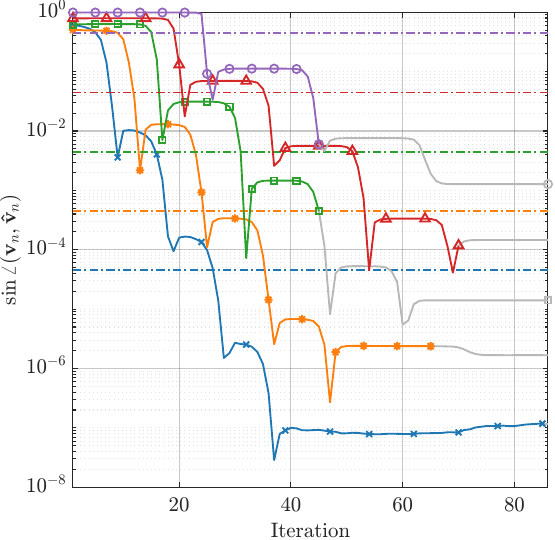}
        \label{subfig:gapSingvec}
    \end{subfigure}
    \caption{Convergence behavior of singular values (\emph{left}) and singular vectors (\emph{right}) of RLOBPCG on matrices with different singular value gaps.
    Dashed lines show predictions from matrix perturbation theory.}
    \label{fig:gapDependency}
\end{figure}

\Cref{fig:gapDependency} displays the convergence histories for both the minimum singular value (\emph{left}) and its associated singular vector (\emph{right}).
As the relative spectral gap becomes smaller, the problem becomes increasingly ill-conditioned (difficult), and the speed of convergence is affected.
Nevertheless, RLOBPCG consistently achieves accuracy levels far better than the dashed lines predicted by perturbation theory for a stable algorithm.

When the gap is small, the singular vector error exhibits a ``spiky'' convergence profile where the error repeatedly drops rapidly, rises to a plateau, then drops again, before settling down at a final level of accuracy.
This behavior is characteristic of nearly degenerate singular subspaces and has been observed in earlier empirical studies of iterative methods \cite{WS15}.

\subsubsection{Effect of the condition number}
Next, we explore how the condition number of the matrix affects the convergence of RLOBPCG. 
In this plot, the singular values decay geometrically from $10^{-8}$ to $10^{-10}$ for the smallest $n-1$ singular values, while the largest singular value is adjusted to yield a prescribed value for the condition number $\kappa_2(\mat{A}) = \sigma_1(\mat{A})/\sigma_n(\mat{A})$.The relative spectral gap remains constant across all test cases so that only the condition number varies between the curves.

\begin{figure}[t]
    \centering
    \begin{subfigure}[t]{0.48\textwidth}
        \centering
        \includegraphics[width=\linewidth]{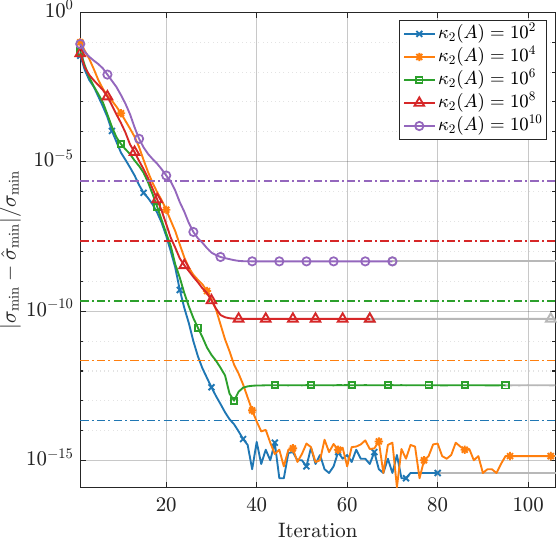}
        \label{subfig:condnoSingval}
    \end{subfigure}%
    \hfill
    \begin{subfigure}[t]{0.48\textwidth}
        \centering
        \includegraphics[width=\linewidth]{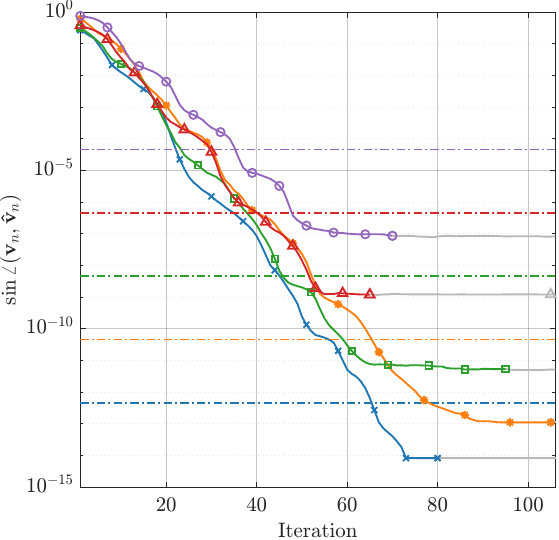}
        \label{subfig:condSingvec}
    \end{subfigure}
    \caption{Convergence behavior of singular values (\emph{left}) and singular vectors (\emph{right}) of RLOBPCG on matrices with different condition numbers.
    Dashed lines show predictions from matrix perturbation theory.}
    \label{fig:condnoDependency}
\end{figure}

\Cref{fig:condnoDependency} shows the results.
As the condition number increases, the maximum attainable accuracy of both quantities deteriorate due to the ill-conditioning of the problem, which amplifies roundoff errors. For the singular values, the limiting accuracy worsens proportionally to the condition number (the observed limiting absolute errors are consistent with the perturbation theory scale $O(u\|\mat{A}\|)$ for singular values, corresponding to relative errors of $O(u\kappa_2(\mat{A}))$ after division by $\sigma_n$).
For the singular vectors, the attained accuracy also worsens with the condition number, though RLOBPCG continues to achieve accuracy levels well below the accuracy level $u/\mathrm{gap}_{\rm abs}$ predicted by perturbation theory. Up to $\kappa_2(\mat A) = 10^{10}$, we saw no case in which the attained accuracy fell short of $u \kappa_2(\mat A)$, although at such condition numbers the stopping criterion may become less reliable.

\subsubsection{Structured problems and the effect of the embedding dimension} \label{subsubsec:sketchsize}
Previous experiments draw sets of singular vectors from the Haar distribution, which is a favorable setting for randomized embeddings, and they use the same embedding dimension. We now vary both. Here we take singular values decaying geometrically from $1$ to $10^{-9}$ over the first $n-1$ values with $\sigma_{\min} = 10^{-10}$, so that $\kappa_2(\mat A) = 10^{10}$ and $\mathrm{gap}_{\rm abs} = 9\times 10^{-10}$. For the first two problems we take $\mat U$ to be the orthonormal factor of $[\mat{G}_1 ; 10^{-7}\mat{G}_2]$ with $\mat{G}_1 \in \real^{n\times n}$ and $\mat{G}_2 \in \real^{(m-n)\times n}$ Gaussian, so that essentially all of the mass of the column space lies on first $n$ of the $m$ rows. This is a coherent column space, on which embeddings built by sampling rows may struggle. We pair it with two choices of right singular vectors:
\begin{itemize}
    \item \textbf{Localized.} $\mat{V} = \Id$, so that $\vec{v}_{\min} = \evec_n$ is supported on a single coordinate.
    \item \textbf{Delocalized.} $\mat{V}$ Haar-random, so that $\vec{v}_{\min}$ is dense.
    \item \textbf{Sparse.} $\mat{A}$ drawn sparse with the same prescribed singular values (MATLAB's \texttt{sprandn}), giving roughly five nonzeros per row.
\end{itemize}
On each problem we run RLOBPCG with embedding dimensions $d/n \in \{1.2, 1.5, 2, 4, 8\}$ at sparsity $\zeta = 4$.

\begin{figure}[!htp]
    \centering
    \begin{subfigure}[t]{0.465\textwidth}
        \centering
        \includegraphics[width=\linewidth]{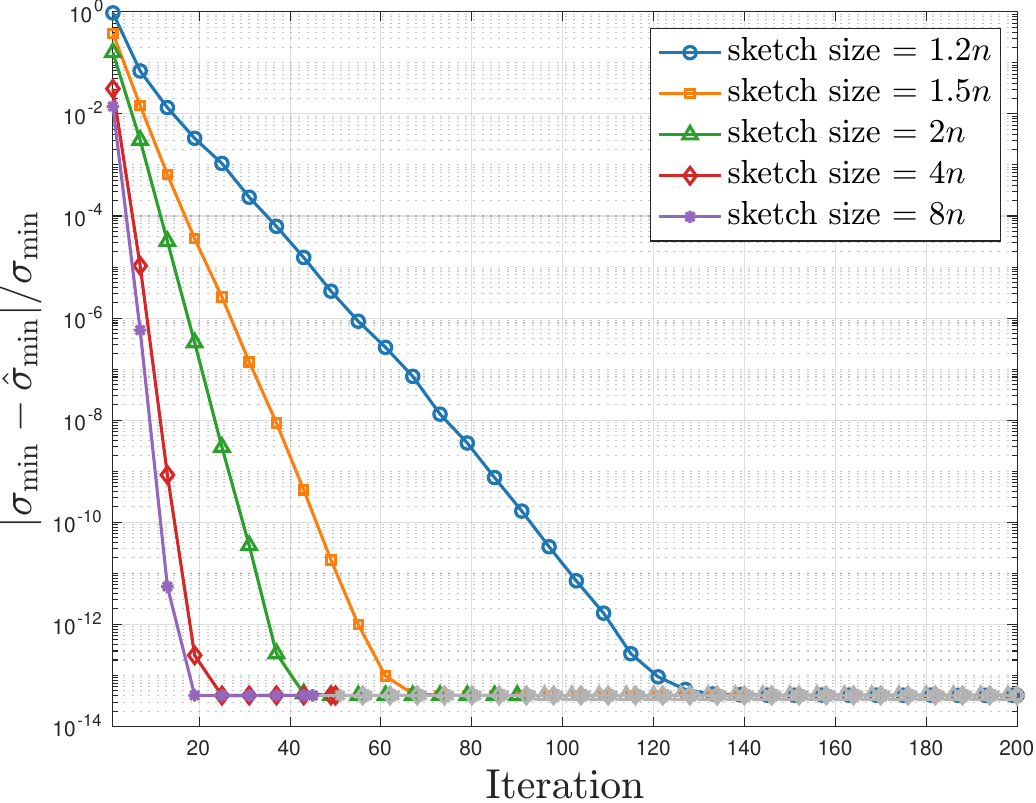}
        \label{subfig:locSingval}
    \end{subfigure}%
    \hfill
    \begin{subfigure}[t]{0.465\textwidth}
        \centering
        \includegraphics[width=\linewidth]{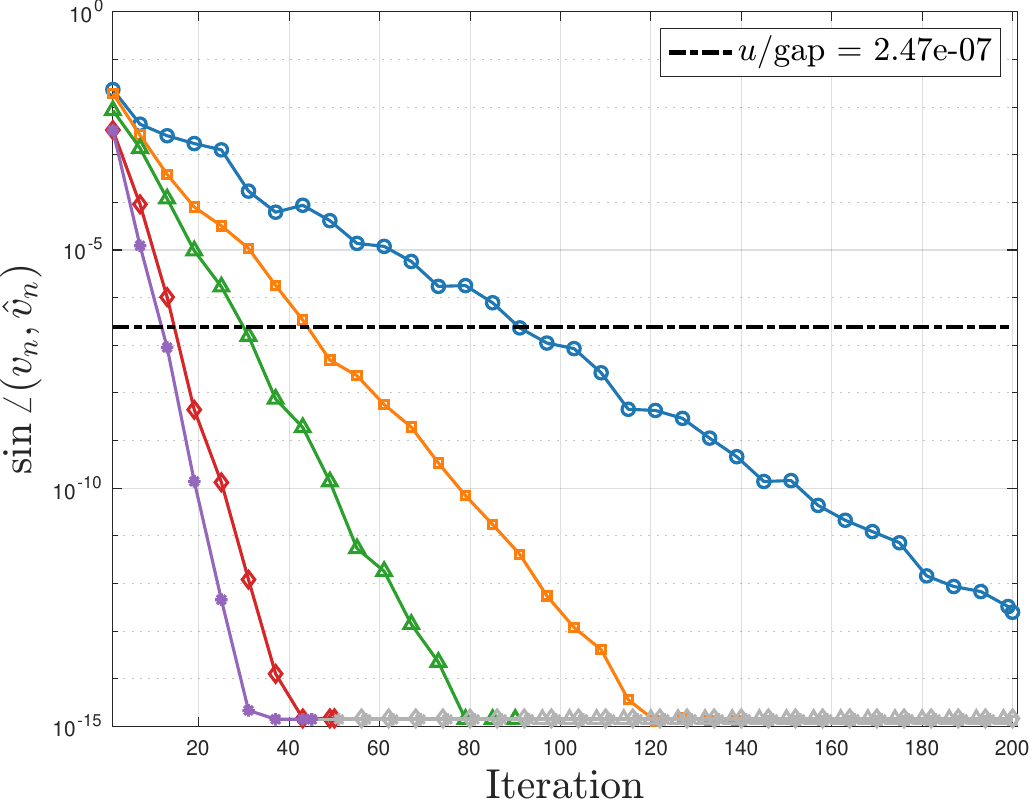}
        \label{subfig:locSingvec}
    \end{subfigure}

    \begin{subfigure}[t]{0.465\textwidth}
        \centering
        \includegraphics[width=\linewidth]{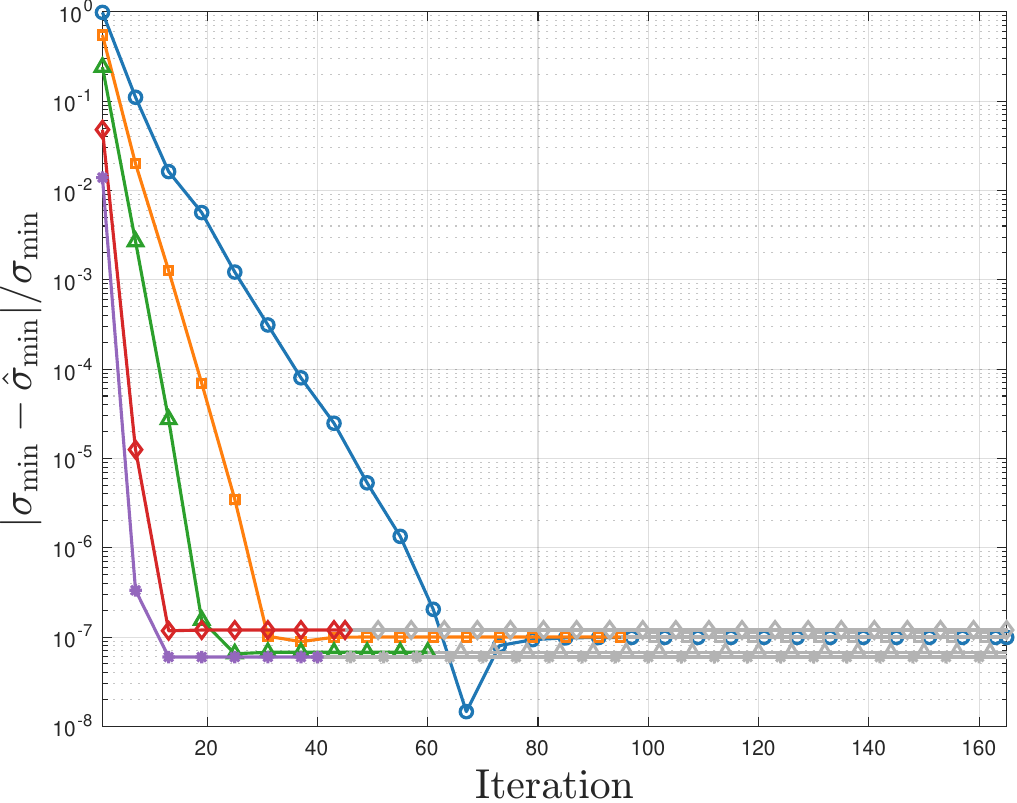}
        \label{subfig:incohSingval}
    \end{subfigure}%
    \hfill
    \begin{subfigure}[t]{0.465\textwidth}
        \centering
        \includegraphics[width=\linewidth]{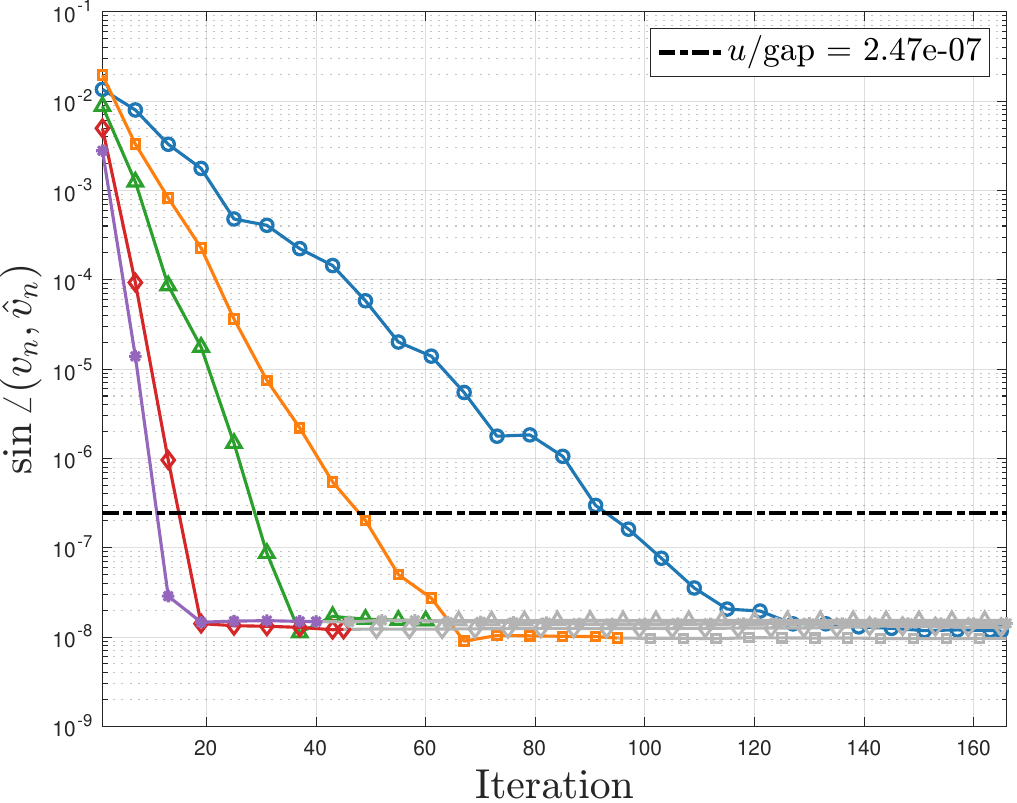}
        \label{subfig:incohSingvec}
    \end{subfigure}

    \begin{subfigure}[t]{0.465\textwidth}
        \centering
        \includegraphics[width=\linewidth]{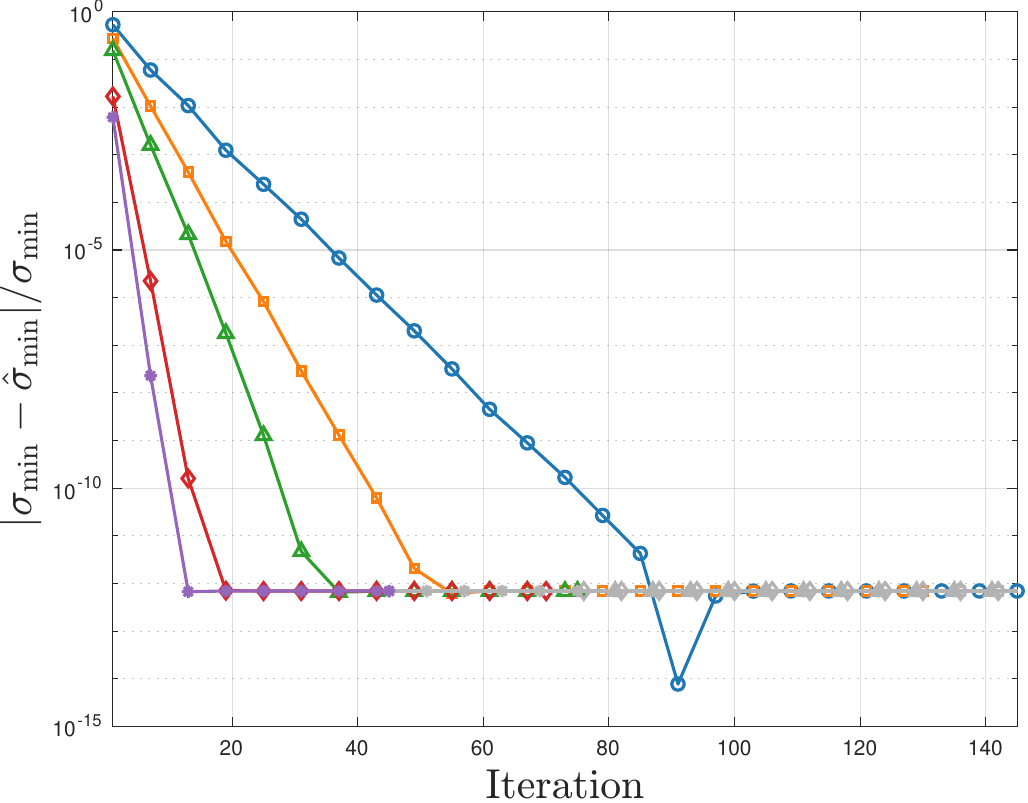}
        \label{subfig:sparseSingval}
    \end{subfigure}%
    \hfill
    \begin{subfigure}[t]{0.465\textwidth}
        \centering
        \includegraphics[width=\linewidth]{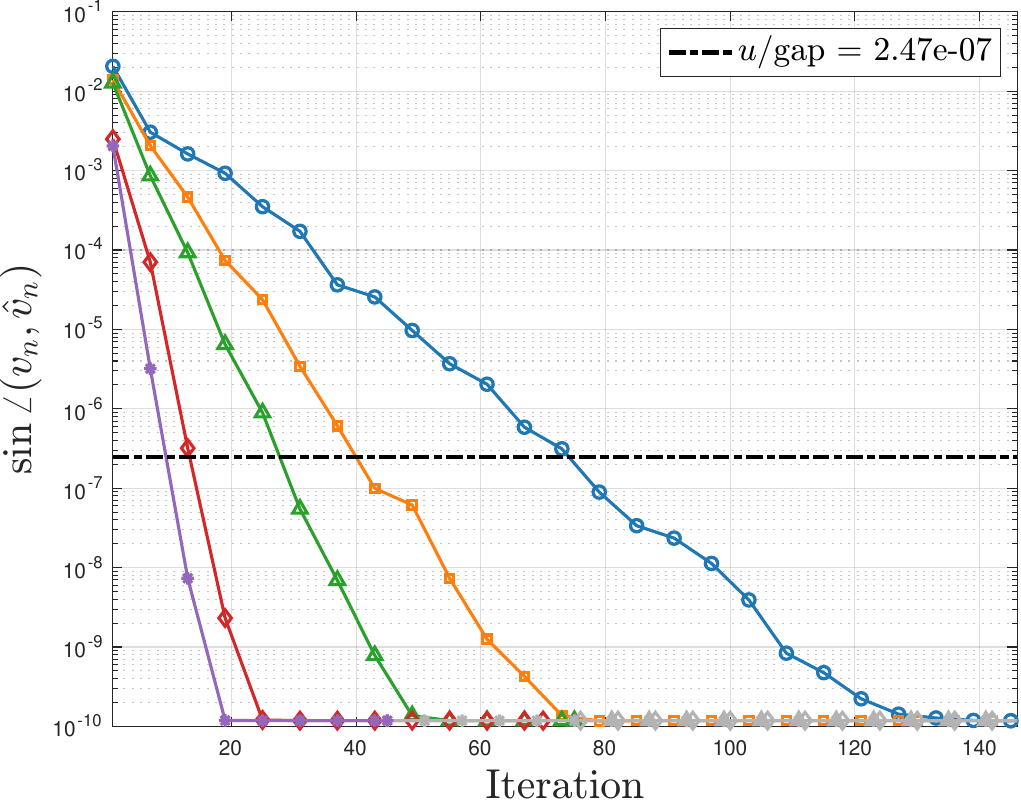}
        \label{subfig:sparseSingvec}
    \end{subfigure}
    \caption{Singular value (\emph{left}) and singular vector (\emph{right}) error for RLOBPCG with a range of embedding dimensions $d$, on problems with a coherent column space and a localized near-nullspace (\emph{top}), a coherent column space and a delocalized near-nullspace (\emph{middle}), and a sparse matrix (\emph{bottom}). The dashed line shows predictions from matrix perturbation theory.}
    \label{fig:structured}
\end{figure}

\Cref{fig:structured} gives two main observations.
First, the examples show no deterioration in convergence attributable to coherence. For a fixed $d$, the pre-asymptotic slopes are similar across the three problems, despite the nearly maximally coherent column space in the first two and the sparsity of the third. This is consistent with the theory. SparseStack is an oblivious embedding, so \cref{fact:sketching_constructions} applies to each fixed matrix without a coherence assumption. In our complexity bounds, sparsity of $\mat{A}$ affects the cost through $\operatorname{nnz}(\mat{A})$, but not the convergence guarantee.

Second, the embedding dimension mainly controls the convergence rate rather than the observed accuracy floor. Increasing $d$ generally improves the quality of the randomized preconditioner and produces faster convergence, as predicted by \cref{eq:thm_convergence_rate}. For each of the three problems, however, all five embedding dimensions eventually reach essentially the same accuracy. Larger $d$ also increases the cost of the SVD of the $d\times n$ sketch in line~\ref{line:sketch-svd}, so these results support $2n\le d\le4n$ as a practical compromise between preprocessing cost and iteration count.

\subsection{Block tests} \label{sec:blocktest}
We next investigate the block implementation of RLOBPCG detailed in \cref{sec:block}.
The block method enables RLOBPCG to compute approximations to multiple singular value/vector triplets at once, and it can improve the reliability and rate of convergence of RLOBPCG \emph{even when only one singular value/vector triplet is required}.
Throughout this section, we will focus on this second use case, requesting only the smallest singular triplet independently of the block size $b$, and evaluate the accuracy of the minimum singular values and vectors computed by block RLOBPCG. We apply the stopping criteria \eqref{eq:est_conv_criterion}--\eqref{eq:stagnation-check} only to the smallest Ritz pair.

\paragraph{Experimental setup}
In this experiment, we vary the block size $b\in\{1,3,5,10,20\}$ to assess how increasing the block dimension influences convergence and computational efficiency.
We generate a test matrix of size $10^6 \times (2\cdot 10^{3})$, with left and right singular vectors drawn from the Haar distribution and \emph{clustered} singular values
\begin{equation*}
    \sigma_i(\mat{A}) = \exp\left(-\log(10^{10})\cdot\left[\frac{i-1}{n-1}\right]^a\right), \quad  i = 1,2,...,n,
\end{equation*}
with $a = 1/256$.
This problem has a small singular value gap $(\sigma_{n-1}^2 - \sigma_n^2) / \sigma_n^2 \approx 9 \times 10^{-5}$, making it a challenging problem for any iterative algorithm.
\paragraph{Measuring cost}
There are several ways of measuring the cost of block LOBPCG.
The simplest method is to count the number of \emph{iterations} as a proxy for the cost of the algorithm.
However, this metric is overly flattering to block methods, as the amount of work performed in each step of the algorithm increases with the block size.
Thus, another natural metric would be the \emph{number of matrix--vector products} with the matrix $\mat{A}$ expended by each algorithm.
Each step of the block method performs a matrix product $\mat{A}\mat{V}$, which may be realized as $b$ matrix--vector products.
But the matrix--vector product count has its own issues, as it ignores the fact that modern computer processors perform matrix--matrix multiplication much faster than iterated matrix--vector multiplication.
Therefore, in addition to these implementation and hardware-insensitive cost metrics, it is appropriate to also report the \emph{wall-clock time} expended by each method.
Given the merits and demerits of measuring cost with each of these metrics, we report all three in this section.

\begin{figure}[!htp]
    \centering
    \begin{subfigure}[t]{0.465\textwidth}
        \centering
        \includegraphics[width=\linewidth]{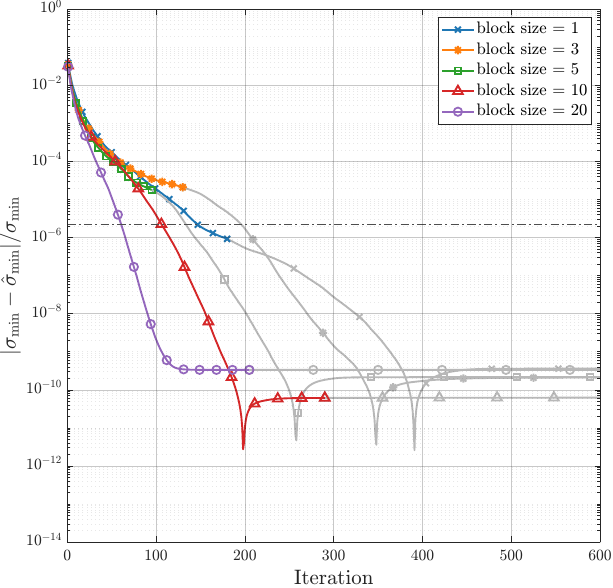}        \label{subfig:clusterIterSingval}
    \end{subfigure}%
    \begin{subfigure}[t]{0.465\textwidth}
        \centering
        \includegraphics[width=\linewidth]{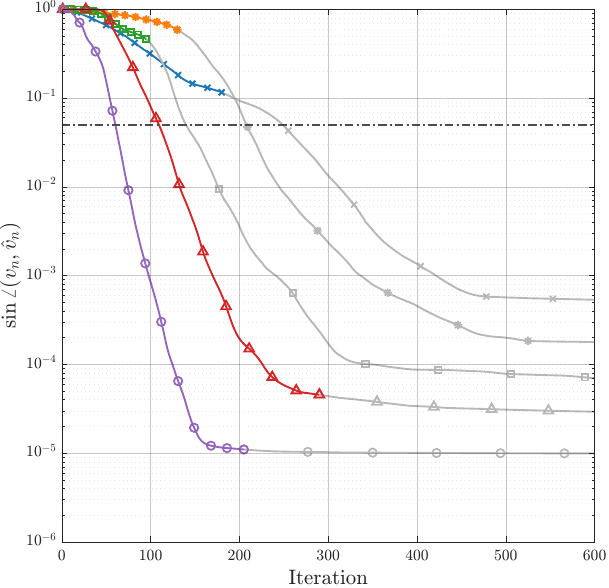}        \label{subfig:clusterIterSingvec}
    \end{subfigure}
    \begin{subfigure}[t]{0.465\textwidth}
        \centering
        \includegraphics[width=\linewidth]{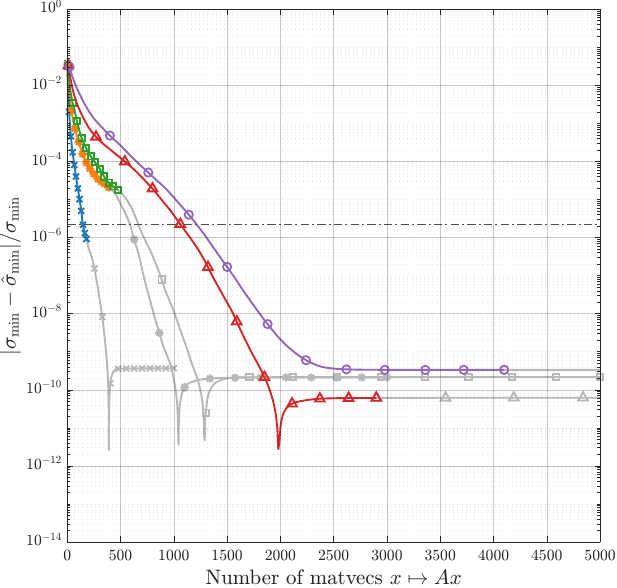}        \label{subfig:clusterMatvecSingval}
    \end{subfigure}%
    \begin{subfigure}[t]{0.465\textwidth}
        \centering
        \includegraphics[width=\linewidth]{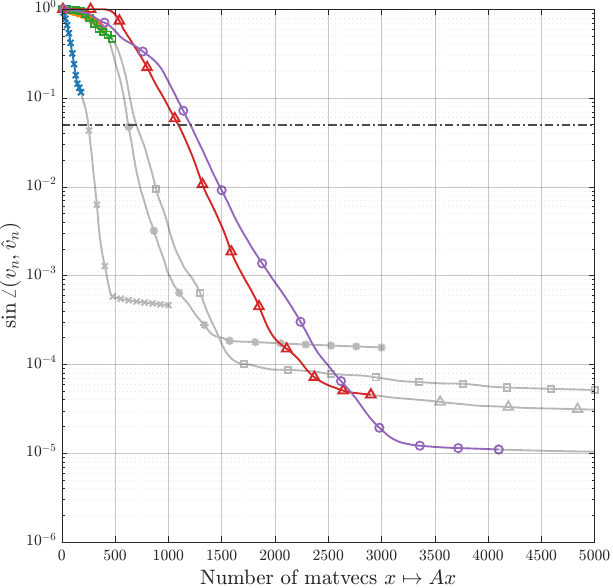}        \label{subfig:clusterMatvecSingvec}
    \end{subfigure}
    \begin{subfigure}[t]{0.465\textwidth}
        \centering
        \includegraphics[width=\linewidth]{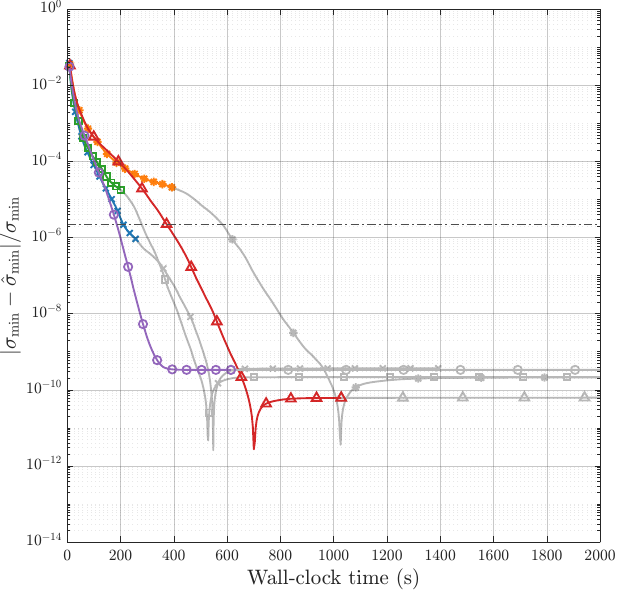}        \label{subfig:clusterTimeSingval}
    \end{subfigure}%
    \begin{subfigure}[t]{0.465\textwidth}
        \centering
        \includegraphics[width=\linewidth]{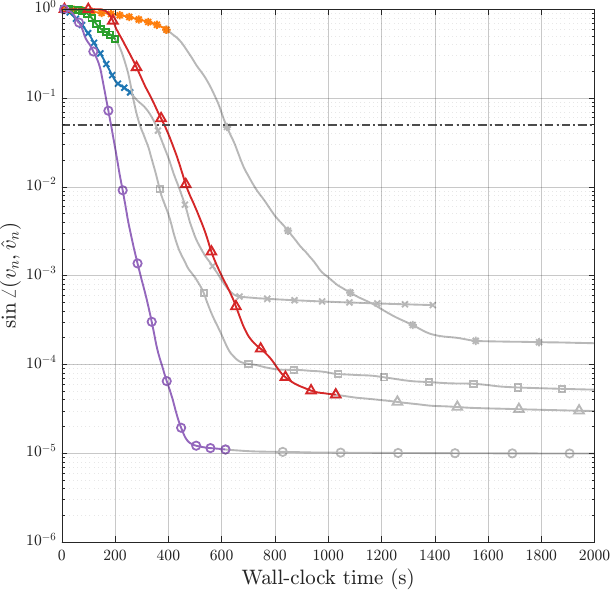}        \label{subfig:clusterTimeSingvec}
    \end{subfigure}
    \caption{Singular value (\emph{left}) and singular vector (\emph{right}) error as a function of iteration count (\emph{top}), total number of matrix--vector products (\emph{middle}), and wall-clock time (\emph{bottom}) for RLOBPCG with various block sizes. Dashed lines show predictions from matrix perturbation theory.}
    \label{fig:clusterRobustness}
\end{figure}

\paragraph{Results}

\Cref{fig:clusterRobustness} presents the results.
We highlight a few conclusions:

\begin{itemize}
    \item \textbf{Which block size is fastest?}
    Larger blocks accelerate convergence per iteration, but require more matrix--vector products. For example, reaching a relative singular-value error of $10^{-8}$ takes 325 iterations for $b = 1$ and 91 iterations for $b = 20$, corresponding to $325$ and $1820$ matrix--vector products, respectively. Wall-clock time gives a different comparison. In this run, $b = 20$ reaches this accuracy fastest taking approximately $275$ seconds, compared with $457$ seconds for $b = 1$ and $417$ seconds for $b = 5$. Thus, the additional work associated with a larger block can be offset by more efficient matrix-matrix operations. The preferred block size depends on the cost metric, the desired accuracy, and the hardware.
    \item \textbf{Stopping criteria.}
    Stopping conditions for RLOBPCG in the presence of both ill-conditioning and small singular value gaps remain challenging, even when the criteria monitor only the requested singular triplet. For $b = 1,3,5$, the gray continuations show substantial improvement after the stopping criteria are first met. For example, at $b =5$ the criteria trigger at iteration $95$ with relative singular-value error approximately $10^{-5}$, whereas an error below $10^{-8}$ is first reached at iteration 202. For $b = 10,20$, the criteria instead trigger much later at a much higher accuracy. Thus, monitoring only the requested triplet does not eliminate premature stopping. Developing stopping criteria that better distinguish slow convergence from stagnation is a natural direction for future work.
    \item \textbf{Attained accuracy.}
    The attained singular-value accuracy is relative insensitive to block size, with errors around $10^{-10}$. Larger blocks improve singular-vector accuracy more substantially. At $600$ iterations, the error for $b = 5$ and $b = 20$ are approximately $1$--$2$ orders of magnitude smaller, respectively, than for $b = 1$.
\end{itemize}
These results show that increasing the block size $b$ can improve convergence and accuracy even when less than $b$ singular triplet is requested. The benefits must be weighed against the additional matrix--vector products, while wall-clock performance can favor a larger block.

\section{Application: Rational approximation} \label{sec:rational}

The AAA-Lawson algorithm~\cite{nakatsukasa2020algorithm} attempts to compute the best (minimax) rational approximation $r_*$ to a function $f:\Omega\rightarrow \mathbb{C}$ in the sense that
\begin{equation*}
    \norm{r - f}_\infty \coloneqq \max_{z \in \Omega} |r(z)-f(z)|
\end{equation*}
is minimized in a domain $\Omega\subseteq \mathbb{C}$ over all type-$(n-1,n-1)$ rational functions (i.e., those that may be expressed as the ratio of two degree-$(n-1)$ polynomials).
Under mild conditions (smoothness of $f$ and $\Omega$), the minimax approximant $r_*$ is that the error $f-r_*$ exhibits a near-circular curve~\cite{trefethen1981rational} with winding number generically equal to $2n+1$; the corresponding effect for real functions is the perfectly equioscillating error curve~\cite[ch.~24]{trefethenatap}. 

Computationally, the AAA-Lawson algorithm proceeds as follows. 
First compute use the AAA algorithm~\cite{nakatsukasa2018aaa} to compute an initial approximant of type $(n-1,n-1)$.
Then refine the approximation $r$ by solving a sequence of the minimum singular vector problem for a matrix of the form $\mat{WA}~\in~\mathbb{C}^{m\times n}$, wherein $\mat{W}~\in~\mathbb{C}^{m\times m}$ is a diagonal ``weight'' matrix that gets updated at each step. This results in an iteratively reweighted nullspace problem; here we employ RLOBPCG for each of these. 

AAA-Lawson provides a natural setting to test RLOBPCG, as the matrix $\mat{A}$ needs to be very tall-skinny to accurately represent the near-circularity of the error curve, and hence to obtain the minimax solution of high accuracy; this is unlike the standard AAA algorithm, for which it is often possible to work with matrices of constant aspect ratio $m=O(n)$, e.g., using an adaptive sampling strategy~\cite{driscoll2024aaa}. 

To evaluate the proposed combination of AAA-Lawson and RLOBPCG, we consider approximating the function 
\begin{equation*}
    f(z) = z\operatorname{sign}(\Re(z)),
\end{equation*}
with $\Omega$ chosen to be the union of unit circles centered at $1.03$ and $-1.03$. 
This problem is challenging because the function $f$ has a singularity along the line $\Re(z) = 0$, which the region $\Omega$ is close to.

We use $10,000$ uniformly distributed points on each circle, resulting in $m=20,000$ points in total, and find a degree $100$ rational approximation. The main cost lies in computing the SVD of a $20,000\times 101$ matrix. 
We compare four algorithms: (i) standard AAA (which does not aim to be minimax), (ii) AAA-Lawson, (iii)  AAA-Lawson where we use sketch-and-solve for the tall-skinny SVDs (which is not expected to yield a minimax approximation), 
(iv) AAA-Lawson where RLOBPCG is used. 
For (ii), (iii), and (iv), we performed $100$ Lawson steps. 


\begin{figure}[htbp]
    \centering
    \begin{subfigure}[t]{0.48\textwidth}
        \centering
        \includegraphics[width=0.98\textwidth]{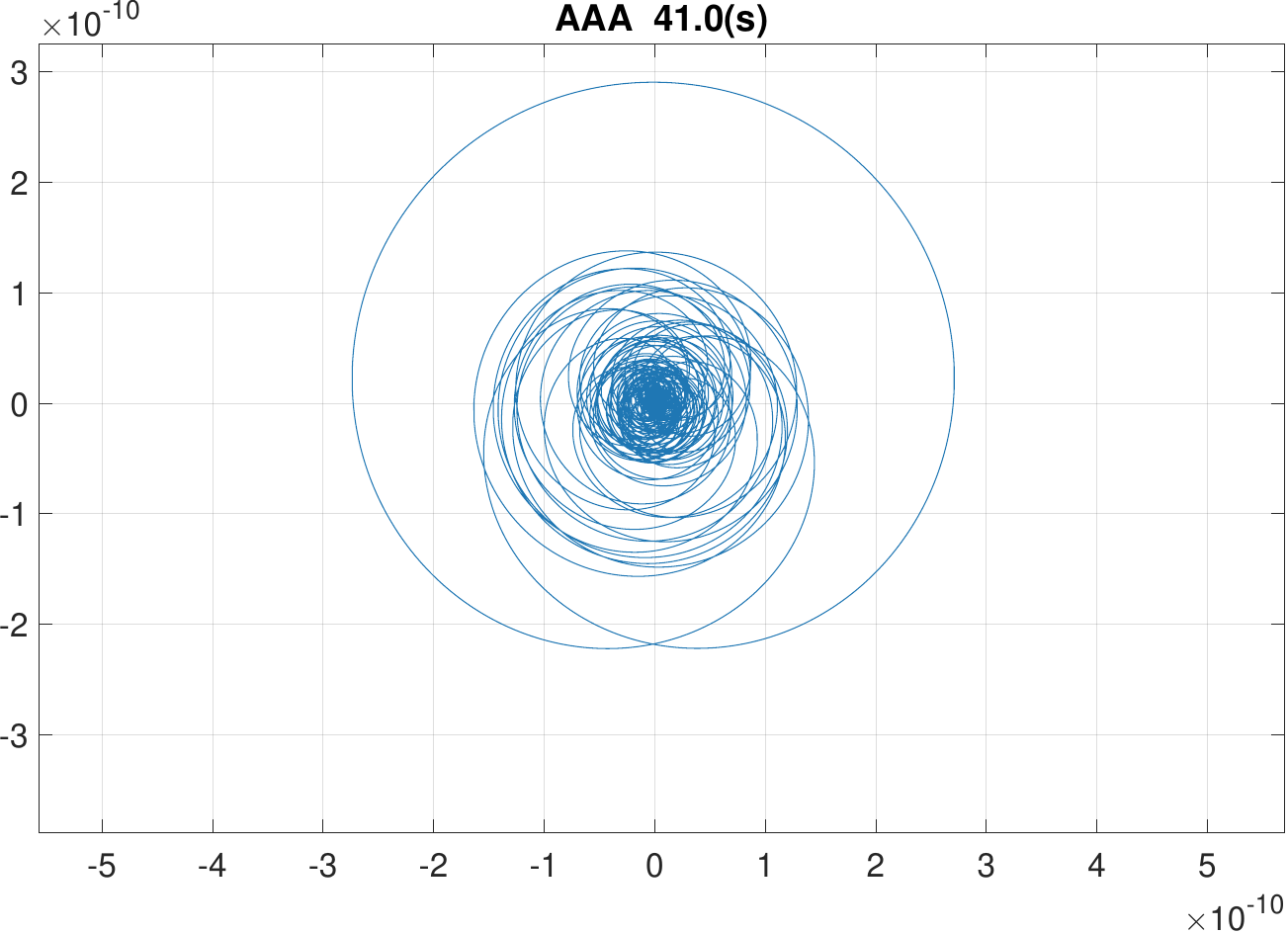}
        
        \label{fig:AAA}
    \end{subfigure}%
    \hfill
    \begin{subfigure}[t]{0.48\textwidth}
        \centering
\includegraphics[width=0.98\textwidth]{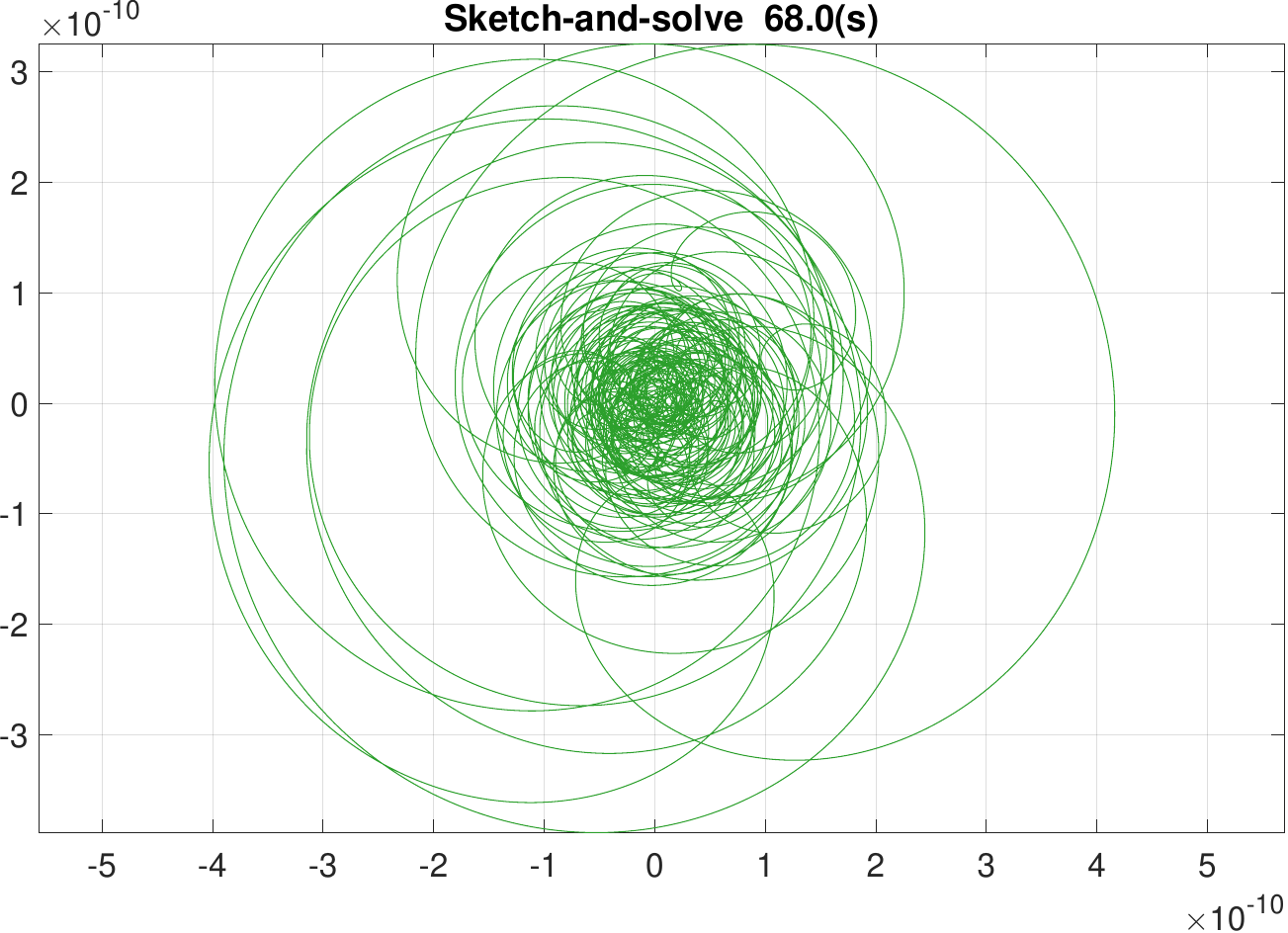}
        
    \end{subfigure}
    \vspace{1em}
    \begin{subfigure}[t]{0.48\textwidth}
        \centering
        \includegraphics[width=0.98\textwidth]{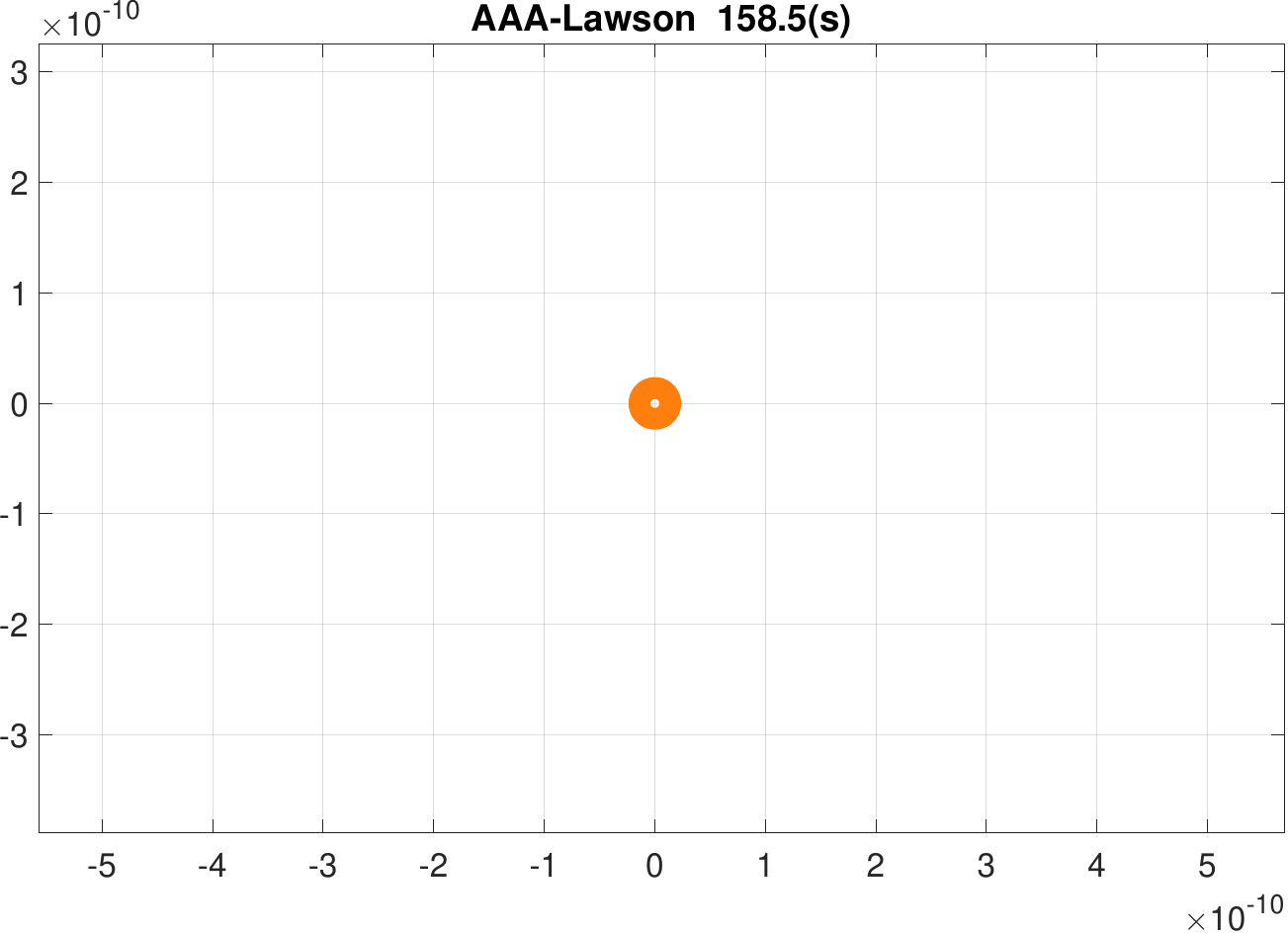}
        \label{fig:AAALaw}
    \end{subfigure}%
    \hfill
    \begin{subfigure}[t]{0.48\textwidth}
        \centering
    \includegraphics[width=0.98\textwidth]{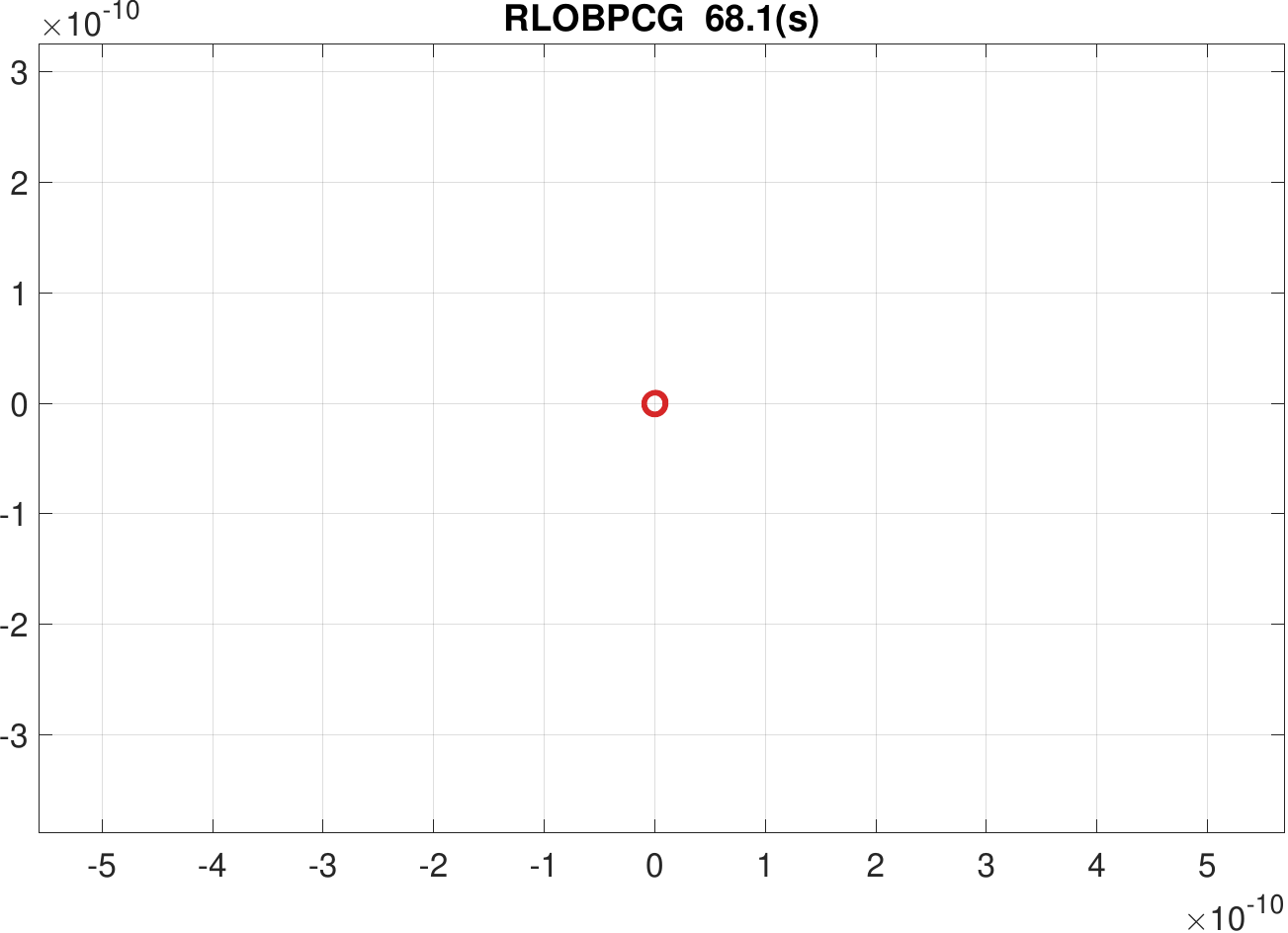}                
        \label{fig:AAALawRLOBPCG}
    \end{subfigure}
    \caption{Error curve $r-f$ of rational approximation to $f(z)=z\operatorname{sign}(\Re(z))$ for AAA (\emph{top left}) and AAA-Lawson with MATLAB's \texttt{svd} (\emph{bottom left}), sketch-and-solve (\emph{top right}), and RLOBPCG (\emph{bottom right}).}
    \label{fig:AAAtotal}
\end{figure}

The resulting execution time and error curves are shown in \cref{fig:AAAtotal}. Observe that the methods based on AAA-Lawson are able to reduce the error by about an order of magnitude as compared with AAA (which is a typical behavior), while the use of RLOBPCG results in more than twice faster runtime compared with standard AAA-Lawson (which uses MATLAB's \texttt{svd}). 
Note that here we used only the stopping criterion~\cref{eq:est_conv_criterion} with $\tau = 10^{-15}$; when we also used the stagnation checks~\cref{eq:stagnation-check}, AAA-Lawson with RLOBPCG further improved the error curve to near-circular (and thus reducing $\|f-r\|_\infty$ by a factor $\approx 2.5$), while the execution time increased to roughly the same as standard AAA-Lawson.

\section{Conclusion}
This work introduces RLOBPCG, a combination of randomized preconditioning with an existing eigensolver for computing a small number of the smallest singular values and associated right singular vectors of very tall matrices. The algorithm combines randomized preconditioning with the LOBPCG eigensolver applied to the normal equations. We prove geometric convergence under the standard subspace embedding assumption and a spectral gap condition. Numerical experiments on the test problems considered demonstrate that randomized preconditioning can significantly improve robustness and efficiency over classical iterative approaches, achieving high accuracy in several ill-conditioned and small-gap regimes. A block variant further improves the method's flexibility.

Several questions remain open. First, the current convergence theory requires a dependence between embedding distortion and the spectral gap; extending the analysis to weaken this dependence would be desirable. Second, while the block implementation performs well empirically, a corresponding convergence theory remains to be developed. Finally, a refined finite-precision analysis could provide additional insight into rounding effects in practice especially when the matrix is ill-conditioned or has small gaps.

\subsection*{Acknowledgments}
The authors thank Joel A.\ Tropp for suggesting to us the idea of using randomization to precondition a tall-skinny nullspace algorithm. 
We also thank Chris Cama\~no and Raphael Meyer for helpful conversations and the referees for helpful and constructive comments.

\bibliographystyle{siamplain}
\bibliography{refs}
	
\end{document}